\documentclass[12pt]{article}
\usepackage{latexsym, amsmath, amsfonts, amscd, amssymb, verbatim, amsxtra,amsthm}
\usepackage{pgf,tikz,pgfplots}
\usepackage{tkz-euclide}
\usepackage{graphicx}
\usepackage{color}
\usepackage{hyperref}
\usepackage[all]{xy}
\usepackage{mathrsfs}
\usepackage{tikz,tkz-tab}
\usepackage{parallel}
\usepackage{array}
\usepackage{xcolor}
\usepackage{verbatim}
\usepackage{tikz-3dplot}
\usepackage{comment}
\usepackage{tikz-cd}

\usetikzlibrary{matrix}
\usetikzlibrary{calc,intersections}
\usetikzlibrary{patterns}
\pgfplotsset{compat=1.15}
\usetikzlibrary{arrows}
\DeclareGraphicsExtensions{.pdf,.png,.jpg}

\begin{document}
	\pagestyle{myheadings}
	\newtheorem{theorem}{Theorem}[section]
	\newtheorem{lemma}[theorem]{Lemma}
	\newtheorem{corollary}[theorem]{Corollary}
	\newtheorem{proposition}[theorem]{Proposition}
	
	\theoremstyle{definition}
	\newtheorem{definition}[theorem]{Definition}
	\newtheorem{example}[theorem]{Example}
	\newtheorem{xca}[theorem]{Exercise}
	\newtheorem{algo}{Algorithm}
	\theoremstyle{remark}
	\newtheorem{remark}[theorem]{Remark}
	
	\numberwithin{equation}{section}
	\renewcommand{\theequation}{\thesection.\arabic{equation}}
	\normalsize
	
	\setcounter{equation}{0}
	
	\title{  \bf  Arrangement of level sets of quadratic constraints and its relation to nonconvex quadratic optimization problems}
	
	\author{Huu-Quang Nguyen\footnote{Department of Mathematics,	Vinh University, Nghe An, Vietnam, email:
			quangdhv@gmail.com.} \  and ~ Ruey-Lin Sheu\footnote{Department of Mathematics, National Cheng Kung University, Tainan, Taiwan email: rsheu@mail.ncku.edu.tw}}
	\maketitle	
	\makeatletter   \renewcommand\@biblabel[1]{#1.}
	\makeatother
	\medskip
	
	\begin{quote}
		\small{\bf Abstract}\ We study a special class of non-convex quadratic programs subject to two (possibly indefinite) quadratic constraints when the level sets of the constraint functions are {\it not} arranged {\it alternatively.} It is shown in the paper that this class of problems admit strong duality following a tight SDP relaxation,	without assuming primal or dual Slater conditions. Our results cover Ye and Zhang's development in 2003 and the generalized trust region subproblems (GTRS) as special cases. Through the novel geometric view and some simple examples, we can explain why the problem becomes very hard when the level sets of the constraints are indeed arranged alternatively.
		
		\medskip
		
		\vspace*{0,05in} {\bf Key words} Non-convex quadratic programming, arrangement of constraint level sets, $\mathcal{S}$-procedure, Strong duality, SDP relaxation, GTRS, CDT problem.

	\end{quote}

	\section{Introduction}\label{sec:intro}
	Consider the following quadratic program
	\begin{equation}\label{key}\begin{array} { c l }  {\rm(P)} & \underset{x \in \mathbb{R}^n}{\min} \, \,    f(x)= x^TA x + 2a^Tx +a_0, \\
	& \text { s.t. }    g(x)=x^TBx + 2b^T x  +b_0\leq 0  \\
	& \hskip 0.85cm    h(x)=x^TC x + 2c^Tx  +c_0\leq 0  \end{array}
	\end{equation} where  $A, B, C$ are symmetric matrices in $\mathbb{R}^{n\times n};$  $a, b, c\in\mathbb{R}^n, $ and $a_0, b_0, c_0\in\mathbb{R}.$ For convenience, let the
	feasible domain be denoted by
	$$\mathcal{D}=\{g\leq 0\}\cap \{h\leq 0\}$$
	where, given an arbitrary function $\phi:\mathbb{R}^n \rightarrow \mathbb{R},$
	$$\{\phi \star \gamma\}=\{x\in\mathbb{R}^n: \phi(x)\star \gamma\},~~\mbox{ for }\star\in\{<, \leq, =, \geq, >\}.$$
	
	It is known that there are no
	general results for (P) but a few special cases do have been studied in literature. They can be classified into several categories based on the setup of the problem data.
	\begin{itemize}
		\item{[CDT]} (The Celis-Dennis-Tapia subproblem)\cite[1984]{CDT}: $B\succ 0$, $C\succeq 0.$
		\item{[HQPD]} ((P) in Homogeneous Quadratic forms with a Positive Definite linear combination) \cite[1998]{Polyak98}:  $a=b=c=0$ and $(\exists \lambda, \beta\in\mathbb{R})$ $\lambda B+\beta C\succ 0.$
		\item{[YZ]} ((P) under conditions assumed in Ye and Zhang \cite[2003]{Ye-Zhang03}): (i) the SDP relaxation of (P), denoted by (SP), and its conic dual (SD) both satisfy the Slater condition (see below for details); (ii) $g(x)$ and $h(x)$ do not have any common root; (iii) $\{h\ge0\} \subset \{g\le0\}.$
		\item{[GTRS]} (The Generalized Trust Region Subproblem)\cite[1993]{More93}\cite[1995]{Stern-Wolkowicz95}\cite[2014]{Pong-Wolkowicz14}:
		$B=-C, b=-c$ with the constraint set $\{l\leq g(x)\leq u\}.$
		\item{[QP1EQC]} ((P) with a single quadratic equality constraint)\cite[1993]{More93}\cite[2016]{Xia-Wang-Sheu16}: $h(x)=-g(x).$  (P) becomes $\min\{f(x)|~g(x)=0\}.$ Note that [QP1EQC] $\subset$ [GTRS].
		\item{[QP1QC]} ((P) with a single quadratic inequality constraint)\cite[1993]{More93}\cite[2007]{Polik-Terlaky07}\cite[2014]{HLS}: $h(x)=-1.$ That is, (P) becomes $\min\{f(x)|~g(x)\le0\}.$
	\end{itemize}
	
	We will use those abbreviations in this paper.
	
	All problems above can be solved in polynomial time. For [CDT], please refer to Bienstock \cite[2016]{Bienstock16}; Sakaue et al. \cite[2016]{Sakaue-Nakatsukasa-Takeda-Iwata16}; and Consolini and Locatelli \cite[2017]{Consolini-Locatelli}. For [HQPD] and [YZ], please see Polyak \cite[1998]{Polyak98} and Ye and Zhang \cite[2003]{Ye-Zhang03}, respectively. For [GTRS]
	and [QP1EQC], they can be resolved by $\mathcal{S}$-lemma with equality. See Wang and Xia \cite[2015]{Wang-Xia} and Xia et al. \cite[2016]{Xia-Wang-Sheu16}.
	For [QP1QC], it can be solved with the classical $\mathcal{S}$-lemma by P$\acute{{\rm o}}$lik and Terlaky \cite[2007]{Polik-Terlaky07} or using a matrix pencil by Hsia et al. \cite[2014]{HLS}.
	
	As far as strong duality is concerned,  [HQPD], [YZ], [QP1QC] admit a tight SDP reformulation under Slater's condition \cite[2020]{X20}; [GTRS], [QP1EQC] requires additionally that $\{g=0\}$ does not separate $\{f<0\}$ \cite[2016]{Xia-Wang-Sheu16}\cite[2018]{Quang-Sheu18}; while [CDT] can do that when and only when the problem data setting does {\it not} satisfy property $\mathcal{I}.$ See Ai and Zhang \cite[2009]{Ai-Zhang09}.
	
	Our main results in this paper add to the above list a new polynomial-solvable subclass in (P), named [Non-Alter], with strong duality.
	\begin{itemize}
		\item{[Non-Alter]} ((P) under Assumptions 1 and 2 below):
	\end{itemize}	
	\noindent {\bf Assumption 1:} $\mathcal{D}\subset
	\{g=0\} \Rightarrow \mathcal{D}=\{g=0\}$ and $\mathcal{D}\subset\{h=0\}\Rightarrow\mathcal{D}=\{h=0\}$.
	
	\noindent {\bf Assumption 2:} Either $\{g=0\}\subset \{h\leq 0\}$ or $\{g=0\}\subset \{h\geq 0\};$ and either $\{h=0\}\subset \{g\leq 0\}$ or $\{h=0\}\subset \{g\geq 0\}.$
	\vskip0.2cm
	Immediately in Section 2, we will show that
	$$\Big({\rm [YZ]}\cup {\rm [GTRS]}  \cup {\rm [QP1EQC]}\cup {\rm [QP1QC]}\Big)\subset \hbox{[Non-Alter]},$$
	but [CDT], [HQPD] may not. We choose particularly the abbreviation [Non-Alter] to stand for ``{\it Non Alternative},'' the geometric meaning of which will be clear in Section 2. The whole idea comes from an earlier paper of us \cite[2018]{Quang-Sheu18}, where we proved that, when a quadratic level set $\{g=0\}$ separates disconnected components of another quadratic hypersurface $\{h<0\}$, the S-lemma with equality \cite[2016]{Xia-Wang-Sheu16} fails to hold.
	
	
	Our major technique to establish strong duality for [Non-Alter] is to develop a new $\mathcal{S}$-procedure involving three general quadratic functions $(f,g,h).$ In Section 3, we first show that if $g,~h$ satisfy Assumptions 1 and 2 along with some constraint qualifications (which can be dealt separately if fail), neither of the following two quadratic systems	has a solution.
	\begin{equation}\label{Unsolvable-system}
	\begin{cases}
	\begin{array}{l}
	g(x)> 0,\\
	h(x)\geq 0;
	\end{array}\end{cases} \hbox{ and }\ \ \  \begin{cases}
	\begin{array}{l}
	g(x)\ge 0,\\
	h(x) >0.
	\end{array}\end{cases}
	\end{equation}
	Then, in Section 4, we show that, under Assumptions 1, 2 and the same constraint qualifications for $(g,h)$, the following two statements are equivalent for any quadratic function $f$ and any given $\gamma\in\mathbb{R}:$
	\begin{itemize}
		\item[](${\rm S_1}$)~~ ($\forall x\in\Bbb R^n$) $~\big(g(x)\le0,~h(x)\le0\big)~\Longrightarrow~ f(x)\ge \gamma.$ 
		\item[](${\rm S_2}$)~~ $(\exists \lambda_1\ge0, \lambda_2\ge0)$ such that
		$f(x)-\gamma+\lambda_1g(x)+\lambda_2h(x)\geq 0,~\forall x\in \mathbb{R}^n.$
	\end{itemize}
	The new $\mathcal{S}$-procedure allows us to compute the optimal value of (P) in the class [Non-Alter]
	through an SDP. We also provide a method to find an optimal solution if (P) is attained.
	%
	%
	
	%

	\section{What are contained in [Non-Alter] and what are not?}
	
	The class [Non-Alter] consists of (P) satisfying Assumptions 1 and 2.
	Assumption 1 relaxes the Slater condition, which we discuss below.
	
	\noindent {\bf $-$ Discussions about Assumption 1:}
	
	Unlike the Slater condition, Assumption 1 allows the feasible domain $\mathcal{D}$ to have an empty interior. When it happens so, either $\mathcal{D}\subset \{g=0\}$ or $\mathcal{D}\subset \{h=0\}$ or both, $\mathcal{D}$ must be the entire $\{g=0\}$ and/or $\{h=0\},$ but cannot be only part of the boundary.
	
	Trivially, [QP1QC] satisfies Assumption 1 due to $\mathcal{D}=\{g\le0\}.$
	In addition, [QP1EQC] also satisfies Assumption 1 since $\mathcal{D}=\{g=0\}=\{h=0\}.$
	
	On the other hand, [QP1QC] satisfies Assumption 2 since we express particularly [QP1QC] with $h(x)=-1.$ Then,
	$\{g=0\}\subset \{h\le 0\}=\mathbb{ R }^n$ and $\{h=0\}=\emptyset\subset \{g\le0\}.$ [QP1EQC] always satisfies Assumption 2 due to $h=-g.$
	
	Therefore,
	$$\Big({\rm [QP1QC]}\cup{\rm [QP1EQC]}\Big) \subset \hbox{[Non-Alter]}.$$
	
	The following example fails Assumption 1. Let $g(x,y,z)=x^2+y^2-1$ and $ h(x,y,z)=-x+1$
	be defined on $\mathbb{R}^3.$ Then, $\mathcal{D}=\{(1,0,t): t \in \mathbb{R}\}=\{h=0\}\subset\{g=0\}$
	but $\mathcal{D}\not=\{g=0\}.$ The example, however, can be solved easily. 
	\vskip0.2cm
	\noindent {\bf $-$ Discussions about Assumption 2:}
	
	Assumption 2 requires that the level set $\{g=0\}$ must lie, either entirely in the sublevel set
	$\{h\le0\},$ or entirely in the superlevel set $\{h\ge0\}.$ Due to symmetry, it also assumes that either $\{h=0\}\subset \{g\leq 0\}$ or $\{h=0\}\subset \{g\geq 0\}.$  Assumption 2 is closely related to a concept for separation between quadratic hypersurfaces, introduced in \cite[2018]{Quang-Sheu18} by the following definition:
	
	\begin{definition}[\cite{Quang-Sheu18}]   \label{DEF1}
		The level set $\{g=0\}$ is said to separate the (sub-, super-) level set $\{h\star \, 0\}$ for $\star\in\{<, \leq, =, \geq, >\}$ if there are non-empty subsets $L^-$ and
		$L^+$ of $\{h \star \, 0\}$  such that $\{h \star \, 0\}=L^-\cup L^+$ and
		$$g(a^-)g(a^+)<0,~\forall~ a^-\in L^-;~\forall~a^+\in
		L^+.$$
	\end{definition}
	
	Let us choose $\star$ to be $=$ for Assumption 2. When $\{h=0\}\subset \{g\leq 0\}$ or $\{h=0\}\subset \{g\geq 0\}$ happens, $\{g=0\}$ cannot separate $\{h=0\}.$ Similarly, $\{g=0\}\subset \{h\leq 0\}$ or $\{g=0\}\subset \{h\geq 0\}$ implies that $\{h=0\}$ cannot separate $\{g=0\}.$ Assumption 2 assumes a symmetry that $\{g=0\}$ and $\{h=0\}$ cannot mutually separate each other.
	
	What types of domain $\mathcal{D}$ satisfy Assumption 2? and what types not?
	
	First, [GTRS] $\subset$ [Non-Alter]. By $l\leq g(x)\leq u,~u\ge l,$ we can write
	$$g_+(x)=g(x)-u,~~g_-(x)=l-g(x).$$
	Then, $$g_+(x)=0~ \Rightarrow ~g(x)=u ~\Rightarrow~ g(x)\ge l ~\Rightarrow ~ g_-(x)\le0.$$
	Similarly, $g_-(x)=0~ \Rightarrow ~g_+(x)\le0.$ So, [GTRS] satisfies Assumption 2.
	As for Assumption 1, if $l\not=u,$ it is automatically satisfied since $\mathcal{D}=\{l\leq g(x)\leq u\}$ consists of an interior point. When $l=u,$ [GTRS] is reduced to [QP1EQC] $\subset$ [Non-Alter].
	
	On the other hand, [CDT] stands in the negative. Take $x\in\mathbb{R}^2$ as an example. With $g(x)=\|x\|^2-1$ and $~h(x)=\|x-0.5\|^2-1,$ the circle $\{h=0\}$ is divided into two parts. One part falls inside $\{g\le0\},$ while the other falling outside $\{g>0\}.$ Then, $\{h=0\}\not\subset\{g\le0\},~\{h=0\}\not\subset\{g\ge0\}.$ Assumption 2 does not hold.
	
	Let us image the arrangement of the level sets for $\{g=0\}$ and $\{h=0\}$ in this case. It appears to be, from left to right in sequence, one piece of $\{g=0\}$, one piece of $\{h=0\}$, one piece of $\{g=0\}$, and finally one piece of $\{h=0\}$. We shall call, informally\footnote{We use the notion of {\it alternative arrangement} of level sets in an informal way because it is complicated to give a formal definition if at all possible. We only intend to capture the important feature of geometry shared by many problems in [Non-Alter] using 2D examples, though the true geometry of a general case could be more sophisticated than that.}, the level sets of $\{g=0\}$ and $\{h=0\}$ are arranged {\it alternatively}. The class [Non-Alter] earns its name by that reason. It excludes (P) having such an alternatively arranged level sets of $\{g=0\}$ and $\{h=0\}.$ {It now suggests that {\it Problems (P) having its feasible domain formed by alternatively arranged $\{g=0\}$ and $\{h=0\}$ is hard, in the sense of practical computation.}
	
	As we have mentioned in Introduction, [CDT] has been proved to be polynomially solvable in computational complexity theory. However, unless the problem data does {\it not} satisfy Property $\mathcal{I}$ \cite[2009]{Ai-Zhang09}
	which implies a convex reformulation for [CDT], there are still no algorithms reported that be able to solve efficiently the type of [CDT] having a positive duality gap against its Lagrange.
	
	There is another reason explaining why (P) outside [Non-Alter] is hard. By the following examples, one observes that, when the level sets of $(g,h)$ are
	arranged alternatively, we are faced with finding a {\it local non-global} minimum
	for [QP1EQC]. In general, identifying a feasible solution (or even a KKT point) to be a local minimum could already be very difficult. See Murty and Kabadi \cite[1987]{Murty-Kabadi}.
	
	\begin{example}	\label{22e}Let	$g(x,y)=-x^2+y^2+9, h(x,y)=1-x.$
		Then, the level set
		$$\{g=0\}=\{g=0\}^-\cup\{g=0\}^+$$
		consists of two branches of a hyperbola, $\{g=0\}^-$ and $\{g=0\}^+$, as depicted below in blue.
		The level set $\{h=0\}$ is a vertical line in red and it separate the two
		branches of $\{g=0\}$ as $h(\{g=0\}^-)h(\{g=0\}^+)<0.$
		Assumption 2 fails for this example.
		
		We also observe the level sets $\{g=0\}$ and $\{h=0\}$ are arranged alternatively: in blue-red-blue sequence (one $g$, one $h$, and $g$ again). We consider
		it as a difficult case to solve if, upon the domain
		$$\mathcal{D}=\{x\in\mathbb{R}^2|~-x^2+y^2+9\le0,~x\ge1\},$$
		the radius $r$ of a ball $f(x,y)=(x+2)^2+(y-8)^2\le r^2$ is to be minimized.
		The solution of (P): $\min \{f(x): x\in \mathcal{D}\}$ attains at point $F$, which is a local non-global solution of (QP1EQC): $\min \{f(x): g(x)=0\},$ whose
		global minimum point $E$ is {\it infeasible} to $\mathcal{D}.$
		
	\end{example}	
	\begin{figure}\begin{center}
			\begin{tikzpicture}[scale=0.31][]
			\clip(-12.291,-9) rectangle (12.31,9);		
			
			\draw [samples=50,domain=-0.99:0.99,rotate around={0.:(0.,0.)},xshift=0.cm,yshift=0.cm,line width=1.pt,color=blue] plot ({3.*(1+(\x)^2)/(1-(\x)^2)},{3.*2*(\x)/(1-(\x)^2)});
			\draw [samples=50,domain=-0.99:0.99,rotate around={0.:(0.,0.)},xshift=0.cm,yshift=0.cm,line width=1.pt,color=blue] plot ({3.*(-1-(\x)^2)/(1-(\x)^2)},{3.*(-2)*(\x)/(1-(\x)^2)});
			
			\draw [line width=1.pt, dash pattern=on 4pt off 4pt, color=red] (1.,-11.16) -- (1.,11);	
			
			\begin{scriptsize}
			
			\draw [fill=black] (-4.440581926585568,3.275407958196115) circle (3.5pt);
			\draw[color=black] (-4.182947397488456,3.9865460658683154) node {$E$};
			\draw [fill=black] (3.383042211047321,1.5913383670979897) circle (3.5pt);
			\draw[color=black] (2.3820145213909356,2.0326883519161187) node {$F$};
			
			\end{scriptsize}
			
			\draw [line width=0.7pt,dash pattern=on 3pt off 3pt] (-2.8,4.42) circle (2.0003999600079982cm);
			\draw [line width=0.7pt,dash pattern=on 3pt off 3pt] (-2.8,4.42) circle (6.7993630302437005cm);

			\draw [->] (-14,0)--(12.3,0);
			\draw [->] (0,-14)--(0,9);

			\draw [line width=0.4pt,domain=14.331030996778363:33.65077316625196] plot(\x,{(--197.62477372416166-0.008044555313814428*\x)/14.094220044445478});
			\draw [line width=0.4pt,domain=13.365780209671895:33.65077316625196] plot(\x,{(--197.73336709554718-0.008659988256844997*\x)/15.172470734920141});
			\draw [line width=0.4pt,domain=12.382395038861212:33.65077316625196] plot(\x,{(--225.778446534663-0.010720603118661742*\x)/18.782708735199268});
			\draw [line width=0.4pt,domain=11.405963117767367:33.65077316625196] plot(\x,{(--225.06536701569473-0.04207095679979922*\x)/20.40877030960836});
			\draw [line width=0.4pt,domain=10.453174825951718:33.65077316625196] plot(\x,{(--204.40788222235426-0.011644394092897414*\x)/20.401208796178366});
			\draw [line width=0.4pt,domain=9.413698864761097:33.65077316625196] plot(\x,{(--167.21810351630705-0.010690006368699656*\x)/18.729102624010558});
			\draw [line width=0.4pt,domain=8.63132732629613:33.65077316625196] plot(\x,{(--162.80965267584295-0.01147510942864649*\x)/20.104618715677745});
			\draw [line width=0.4pt,domain=7.537557995752271:33.65077316625196] plot(\x,{(--149.21502846724823-0.012308988201970017*\x)/21.565590822045124});
			\draw [line width=0.4pt,domain=6.743551755156014:33.65077316625196] plot(\x,{(--137.06430930082593-0.04713429732138508*\x)/22.6420378998789});
			\draw [line width=0.4pt,domain=5.828109815886611:33.65077316625196] plot(\x,{(--117.08154676398951-0.013365290891462323*\x)/23.416254029454812});
			\draw [line width=0.4pt,domain=5.029219296594739:33.65077316625196] plot(\x,{(--97.92724441768769-0.013837404572908163*\x)/24.24340653853867});
			\draw [line width=0.4pt,domain=4.235180822591189:33.65077316625196] plot(\x,{(--74.82392382092894-0.0142744857979058*\x)/25.009181490926462});
			\draw [line width=0.4pt,domain=3.597709066605339:33.65077316625196] plot(\x,{(--51.768635118061056-0.014864056688303817*\x)/26.04212135359295});
			\draw [line width=0.4pt,domain=3.1704238336514283:33.65077316625196] plot(\x,{(--26.814230352935418-0.014898343052654806*\x)/26.10219174217714});
			\draw [line width=0.4pt,domain=3.0000007774873927:33.65077316625196] plot(\x,{(--0.1061067565727547-0.015640529658010497*\x)/27.402517356441624});
			
			\draw [line width=0.4pt] (3.157210052089254,-0.9838573641607953)-- (26.723765173277265,-1.055585752929745);
			\draw [line width=0.4pt] (3.551664268775887,-1.9011362597403048)-- (27.12459836243828,-2.0242659600688686);
			\draw [line width=0.4pt] (4.223922062297646,-2.97346894861285)-- (27.358417722782207,-2.959543401444574);
			\draw [line width=0.4pt] (4.927419709772105,-3.908895623606573)-- (27.558834317362717,-3.9282236085836972);
			\draw [line width=0.4pt] (5.76098237598679,-4.91822304663284)-- (27.759250911943223,-4.896903815722821);
			\draw [line width=0.4pt] (6.635191327600963,-5.918256834057731)-- (28.39390346144817,-5.898986788625362);
			\draw [line width=0.4pt] (7.5852149442153,-6.9667414011105)-- (28.427306227211584,-6.96787529305474);
			\draw [line width=0.4pt] (8.436288814757622,-7.884856940110238)-- (28.72793111908235,-7.903152734430446);
			\draw [line width=0.4pt] (9.506752041498748,-9.020994090372781)-- (29.36258366858729,-9.038846770386659);
			\draw [line width=0.4pt] (10.33576017503593,-9.890800695386535)-- (29.663208560458052,-9.87391591447211);
			\draw [line width=0.4pt] (11.20390096969852,-10.7947856365382)-- (29.663208560458052,-10.842596121611233);
			\draw [line width=0.4pt] (12.238374269921497,-11.86498229120788)-- (29.663208560458052,-11.811276328750358);
			\draw [line width=0.4pt] (13.177281262189322,-12.831240838782737)-- (29.930430686565398,-12.846762067416318);
			\draw [line width=0.4pt] (14.314916986009768,-13.99702998197657)-- (30.832305362177685,-13.949053337609113);
			\draw [line width=0.4pt] (15.20794099653008,-14.909106926772653)-- (30.93251365946794,-14.884330778984818);
			\draw [line width=0.4pt] (16.233473502946836,-15.953860409658759)-- (30.865708127941105,-15.919816517650778);
			\draw [line width=0.4pt] (17.16159020841263,-16.897342349656178)-- (31.166333019811866,-16.855093959026483);
			\draw [line width=0.4pt] (18.272276030695927,-18.024318887046604)-- (31.366749614392376,-17.92398246345586);
			\draw [line width=0.4pt] (19.201461947164855,-18.96565688048848)-- (31.767582803553392,-18.9260654363584);
			

			\pgftransformrotate{45}
			\pgftransformshift{\pgfpoint{10.0cm}{0.2cm}}
			\pgfnode{rectangle}{center}
			{\color{blue} \Huge{$\{g=0\}^+$}}{}{}
			
			\pgftransformrotate{-315}
			\pgftransformshift{\pgfpoint{-12cm}{4.9cm}}
			\pgfnode{rectangle}{center}
			{\color{red} \Huge{$\{h=0\}$}}{}{}
			
			\pgftransformrotate{221}
			\pgftransformshift{\pgfpoint{-14.5cm}{1.2cm}}
			\pgfnode{rectangle}{center}
			{\color{blue} \Huge{$\{g=0\}^-$}}{}{}

			\pgftransformrotate{45}
			\pgftransformshift{\pgfpoint{4.5cm}{0.5cm}}
			\pgfnode{rectangle}{center}
			{\color{black} \Huge{$\{f=8.6\}$}}{}{}

			\pgftransformrotate{5}
			\pgftransformshift{\pgfpoint{-3.9cm}{-6.3cm}}
			\pgfnode{rectangle}{center}
			{\color{black} \Huge{$\{f=19.5\}$}}{}{}


			\end{tikzpicture}
		\end{center}	
		\caption{Graphic representation for Example \ref{22e}}
		\label{xxsd}
	\end{figure}
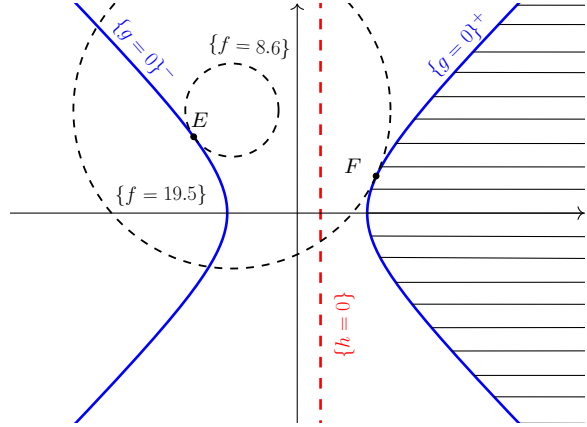
	
	
	\begin{example} \label{23e}
		Let $g(x,y)=-x^2+3y^2+1$ (blue) and $h(x,y)=x^2-3y^2-2x-y-1$ (red).
		Both $\{g=0\}$ and $\{h=0\}$ have two branches:
		$$\{g=0\}=\{g=0\}^-\cup\{g=0\}^+;~~\{h=0\}=\{h=0\}^-\cup\{h=0\}^+.$$
		We can see that they mutually separate each other:
		$$g(\{h=0\}^-)g(\{h=0\}^+)<0;~~h(\{g=0\}^-)h(\{g=0\}^+)<0.$$
		Thus, $\{h=0\}\not\subset\{g\le0\}$ and $\{h=0\}\not\subset\{g\ge0\}.$ The example does not belong to [Non-Alter].
		We observe that the level sets of $\{g=0\}$ and $\{h=0\}$ (4 pieces in total) are arranged alternatively: blue-red-blue-red.
		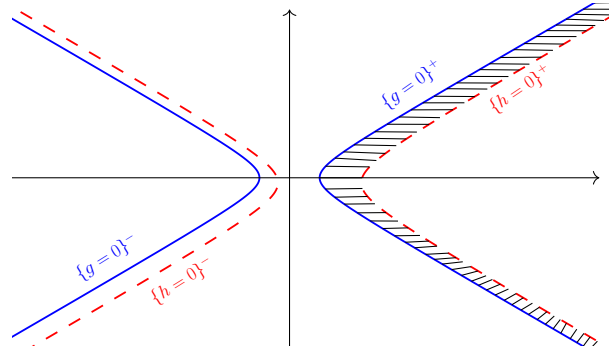
\begin{figure}	\begin{center}
				\begin{tikzpicture}[scale=0.4]
				\clip(-9.210374688934625,-5.593335307413431) rectangle (10.607608492016908,5.796310198880546);
				\draw [samples=50,domain=-0.99:0.99,rotate around={0.:(0.,0.)},xshift=0.cm,yshift=0.cm,line width=0.7pt,color=blue] plot ({1.*(1+(\x)^2)/(1-(\x)^2)},{0.5773502691896257*2*(\x)/(1-(\x)^2)});
				\draw [samples=50,domain=-0.99:0.99,rotate around={0.:(0.,0.)},xshift=0.cm,yshift=0.cm,line width=0.7pt,color=blue] plot ({1.*(-1-(\x)^2)/(1-(\x)^2)},{0.5773502691896257*(-2)*(\x)/(1-(\x)^2)});
				\draw [samples=50,domain=-0.99:0.99,rotate around={0.:(1.,-0.16666666666666666)},xshift=1.cm,yshift=-0.16666666666666666cm,line width=0.7pt,color=red,dash pattern=on 5pt off 5pt] plot ({1.3844373104863459*(1+(\x)^2)/(1-(\x)^2)},{0.7993052538854533*2*(\x)/(1-(\x)^2)});
				\draw [samples=50,domain=-0.99:0.99,rotate around={0.:(1.,-0.16666666666666666)},xshift=1.cm,yshift=-0.16666666666666666cm,line width=0.7pt,color=red,dash pattern=on 5pt off 5pt] plot ({1.3844373104863459*(-1-(\x)^2)/(1-(\x)^2)},{0.7993052538854533*(-2)*(\x)/(1-(\x)^2)});


				\draw [line width=0.4pt] (10.122782121607951,5.815803695722704)-- (11.48677267761534,5.83488131490478);
				\draw [line width=0.4pt] (9.763814088179458,5.607496932634967)-- (11.0918484138993,5.604778494424709);
				\draw [line width=0.4pt] (9.470481843414627,5.437218233203634)-- (10.766951030449748,5.415348076074427);
				\draw [line width=0.4pt] (9.189537320351448,5.274074837104358)-- (10.450956641382144,5.230984512955156);
				\draw [line width=0.4pt] (8.855403928417504,5.079966493834125)-- (10.162435021032492,5.0625317692944956);
				\draw [line width=0.4pt] (8.537415433681117,4.89514937760693)-- (9.86692310664568,4.889868995638781);
				\draw [line width=0.4pt] (8.197192395915401,4.697302884835556)-- (9.536860679587642,4.696847960843734);
				\draw [line width=0.4pt] (7.795936078689148,4.463803286548632)-- (9.110659370900814,4.447302295955228);
				\draw [line width=0.4pt] (7.447925556875667,4.261126419964628)-- (8.754448772473776,4.238436778893318);
				\draw [line width=0.4pt] (7.,4.)-- (8.367313717393841,4.011077730486449);
				\draw [line width=0.4pt] (6.6441715373566925,3.792317471670886)-- (7.973516495823516,3.77935519722174);
				\draw [line width=0.4pt] (6.3114694488744725,3.5978996930651936)-- (7.643359011290526,3.5846687527197143);
				\draw [line width=0.4pt] (5.9232654539408305,3.370710787447699)-- (7.275519424730242,3.36723945279472);
				\draw [line width=0.4pt] (5.557197305081358,3.156075721926598)-- (6.9104457134637265,3.15079731782184);
				\draw [line width=0.4pt] (5.228590772422477,2.963003063642329)-- (6.574343995492924,2.9508453923007267);
				\draw [line width=0.4pt] (4.801411699264947,2.711306837783723)-- (6.15736541522897,2.7016517349717812);
				\draw [line width=0.4pt] (4.431740198240856,2.4926773547535794)-- (5.812137671205244,2.4941604165266997);
				\draw [line width=0.4pt] (4.03716757911524,2.2582236427377467)-- (5.431827586648835,2.264000531937794);
				\draw [line width=0.4pt] (3.6433349340118175,2.0226953503836484)-- (5.,2.);
				\draw [line width=0.4pt] (3.242025993909698,1.7805179532171873)-- (4.639818779264541,1.7768365345330073);
				\draw [line width=0.4pt] (2.833961595879445,1.5309407050263322)-- (4.232755834307586,1.5199516084213929);
				\draw [line width=0.4pt] (2.485788715097846,1.313918507887482)-- (3.8816087545939704,1.2924429188781095);
				\draw [line width=0.4pt] (2.1436836463448854,1.0947419141220003)-- (3.5478695639737055,1.0682385453984662);
				\draw [line width=0.4pt] (1.7965127521056203,0.8616762865636007)-- (3.1855682493145947,0.8097283875583801);
				\draw [line width=0.4pt] (1.509729401038478,0.6530142582312684)-- (2.887008829487396,0.5736340956132606);
				\draw [line width=0.4pt] (1.2001509341995442,0.38312846620945956)-- (2.6600254753949164,0.36217417652828787);
				\draw [line width=0.4pt] (1.,0.)-- (2.414213562373095,0.);
				\draw [line width=0.4pt] (1.1912060362297066,-0.3737074706014891)-- (2.4228745961255242,-0.35632025405931933);
				\draw [line width=0.4pt] (1.5011001025706272,-0.646349110501321)-- (2.594104461354883,-0.6229203852067184);
				\draw [line width=0.4pt] (1.8169666980545305,-0.8758553879568597)-- (2.8378006418964015,-0.864483388709582);
				\draw [line width=0.4pt] (2.254778206454758,-1.1667654377670682)-- (3.206383314990094,-1.158543787080281);
				\draw [line width=0.4pt] (2.6798666702259752,-1.4354657049413213)-- (3.604586969291314,-1.440402198339215);
				\draw [line width=0.4pt] (3.089981447236544,-1.6879953933849117)-- (3.943178332141,-1.6661812579169608);
				\draw [line width=0.4pt] (3.488925495437887,-1.9298187404275828)-- (4.335287914616012,-1.9185683768773125);
				\draw [line width=0.4pt] (4.6329865877898,-2.105904049453511)-- (3.8469639789277026,-2.1446935954557906);
				\draw [line width=0.4pt] (4.9048652652523135,-2.2746914983600175)-- (4.155443340914391,-2.328641188873156);
				\draw [line width=0.4pt] (5.2834713717523645,-2.506998575175616)-- (4.522630886635373,-2.5465133114874448);
				\draw [line width=0.4pt] (5.5320103481558185,-2.658148898132156)-- (4.806556214069172,-2.714343544886865);
				\draw [line width=0.4pt] (5.757720525548757,-2.794672116468301)-- (5.071471028681786,-2.8705294746648327);
				\draw [line width=0.4pt] (6.0262371570756805,-2.9563137944411664)-- (5.383196809063007,-3.0538941198141396);
				\draw [line width=0.4pt] (6.266285466340304,-3.1002137770389337)-- (5.711941914103545,-3.2468590170652036);
				\draw [line width=0.4pt] (6.532577972619058,-3.259278958157629)-- (6.042681238128347,-3.440639308114522);
				\draw [line width=0.4pt] (6.742319154125613,-3.3842000180178338)-- (6.276250105975964,-3.577304170869609);
				\draw [line width=0.4pt] (7.0047932372709285,-3.540135462812458)-- (6.542526178265901,-3.732945610519688);
				\draw [line width=0.4pt] (7.275641100677519,-3.700644810243971)-- (6.83106306812119,-3.9014280565219224);
				\draw [line width=0.4pt] (7.532211669562931,-3.852365272056037)-- (7.132929077411756,-4.077526914825787);
				\draw [line width=0.4pt] (7.769710269265239,-3.9925569971705164)-- (7.386161844413441,-4.225138529952557);
				\draw [line width=0.4pt] (8.072586889408639,-4.171031728344913)-- (7.689908541012141,-4.4020712310279375);
				\draw [line width=0.4pt] (8.304025950474756,-4.307203202694416)-- (7.938779933095524,-4.546948677450999);
				\draw [line width=0.4pt] (8.554851915035592,-4.454600154624483)-- (8.230229418783555,-4.716519772240626);
				\draw [line width=0.4pt] (8.798863889831871,-4.597829467612443)-- (8.500440875063322,-4.873653491664726);
				\draw [line width=0.4pt] (8.998705230585404,-4.715022365859054)-- (8.760039061518313,-5.024549212261952);
				\draw [line width=0.4pt] (9.202960375463734,-4.834710234542523)-- (9.015199401506132,-5.172807756235697);
				\draw [line width=0.4pt] (9.54692018840182,-5.036066973352079)-- (9.32807316232162,-5.354529201888338);
				\draw [line width=0.4pt] (9.78490303512341,-5.175254447330787)-- (9.549613023508577,-5.483159335600983);
				\draw [line width=0.4pt] (9.98928801627844,-5.294714780486079)-- (9.758801695790591,-5.604587720126116);
				\draw [line width=0.4pt] (10.231746261098763,-5.436343202047952)-- (9.961823067229476,-5.722409714528348);

				\draw [->] (-16,0)--(10.3,0);
				\draw [->] (0,-7.6)--(0,5.6);
				
				\pgftransformrotate{35}
				\pgftransformshift{\pgfpoint{5.1cm}{0.2cm}}
				\pgfnode{rectangle}{center}
				{\color{blue} \Large{$\{g=0\}^+$}}{}{}
				
				\pgftransformrotate{-1}
				\pgftransformshift{\pgfpoint{2.9cm}{-2.151cm}}
				\pgfnode{rectangle}{center}
				{\color{red} \Large{$\{h=0\}^+$}}{}{}

				\pgftransformrotate{-362}
				\pgftransformshift{\pgfpoint{-12.9cm}{0.551cm}}
				\pgfnode{rectangle}{center}
				{\color{red} \Large{$\{h=0\}^-$}}{}{}

				\pgftransformrotate{0}
				\pgftransformshift{\pgfpoint{-1.7cm}{2.0cm}}
				\pgfnode{rectangle}{center}
				{\color{blue} \Large{$\{g=0\}^-$}}{}{}

				\label{pic6}
				\end{tikzpicture}
			\end{center}
			
			\caption{Graphic representation for Example \ref{23e}}
		\end{figure}
	\end{example}

	
	\begin{example}\label{24e}	Let $g(x,y)=-x^2+y^2+9$ (blue) and $h(x,y)=-(-x^2+y^2+49)$ (red). Since $\{g=0\}\subset\{h\le0\}$ and $\{h=0\}\subset\{g\le0\},$ Assumption 2 holds. The level sets $\{g=0\},~\{h=0\}$ are arranged in the manner: red-blue-blue-red, which is not alternative.
		
		Consider $f_1(x,y)=(x+2)^2+(y-8)^2$, $f_2(x,y)=(x-14)^2+(y-4)^2$ as objective functions. The solution to $\min \{f_1(x): x\in\mathcal{D} \}$ is at $G$, whereas
		the solution to $\min \{f_2(x): x\in\mathcal{D} \}$ is at $H$. The point $G$ is also the global minimum for $\min \{f_1(x): g=0 \}$ and the point $H$ is the global minimum for $\min \{f_1(x): h=0 \}.$ This example shows that [Non-Alter] can be solved in polynomial time as it reduces to [QP1EQC].
		
	\end{example}	
	\begin{figure}\begin{center}
			\begin{tikzpicture}[scale=0.27][]
			\clip(-16,-8) rectangle (19.2,14);
			
			\draw [->] (-18,0)--(19.2,0);
			\draw [->] (0,-18)--(0,14);
			
			\draw [line width=0.6pt,dash pattern=on 4pt off 4pt] (-2.,8.) circle (4.86867cm);
			\draw [line width=0.6pt,dash pattern=on 4pt off 4pt,color=black] (14.,4.) circle (5.08cm);
			
			\draw [samples=50,domain=-0.99:0.99,rotate around={0.:(0.,0.)},xshift=0.cm,yshift=0.cm,line width=1.pt,color=blue] plot ({3.*(1+(\x)^2)/(1-(\x)^2)},{3.*2*(\x)/(1-(\x)^2)});
			
			\draw [samples=50,domain=-0.99:0.99,rotate around={0.:(0.,0.)},xshift=0.cm,yshift=0.cm,line width=1.pt,color=blue] plot ({3.*(-1-(\x)^2)/(1-(\x)^2)},{3.*(-2)*(\x)/(1-(\x)^2)});
			
			\draw [samples=50,domain=-0.9:0.9,rotate around={0.:(0.,0.)},xshift=0.cm,yshift=0.cm,,color=red,line width=1.pt,dash pattern=on 5pt off 5pt] plot ({7.*(1+(\x)^2)/(1-(\x)^2)},{7.*2*(\x)/(1-(\x)^2)});
			
			\draw [samples=50,domain=-0.9:0.9,rotate around={0.:(0.,0.)},xshift=0.cm,yshift=0.cm,,color=red,line width=1.pt,dash pattern=on 5pt off 5pt] plot ({7.*(-1-(\x)^2)/(1-(\x)^2)},{7.*(-2)*(\x)/(1-(\x)^2)});
			
			\begin{scriptsize}
			
			\draw [fill=black] (-5.8675,5.0426) circle (3.5pt);
			\draw[color=blue] (-5.2,5.81598) node {$G$};
			
			\draw [fill=black] (9.89,6.987) circle (3.5pt);
			\draw[color=red] (10.8, 6.5917) node {$H$};
			
			\draw [fill=black] (-7,0) circle (3.5pt);
			\draw[color=black] (-6.6,0.3) node {$-7$};
			
			\draw [fill=black] (7, 0) circle (3.5pt);
			\draw[color=black] (7.4, 0.3) node {$7$};
			
			\draw [fill=black] (-3,0) circle (3.5pt);
			\draw[color=black] (-2.6,0.3) node {$-3$};
			
			\draw [fill=black] (3, 0) circle (3.5pt);
			\draw[color=black] (3.4, 0.3) node {$3$};
			
			\end{scriptsize}

			\draw [line width=0.5pt] (3.3,1.3747727084867516)-- (3.3,-1.3747727084867516);
			\draw [line width=0.5pt] (3.8,2.33238075793812)-- (3.8,-2.3323807579381204);
			\draw [line width=0.5pt] (4.3,3.080584360149872)-- (4.3,-3.0805843601498726);
			\draw [line width=0.5pt] (4.8,3.746998799039039)-- (4.8,-3.7469987990390385);
			\draw [line width=0.5pt] (5.3,4.369210454990696)-- (5.3,-4.369210454990696);
			\draw [line width=0.5pt] (5.8,4.963869458396343)-- (5.8,-4.963869458396342);
			\draw [line width=0.5pt] (6.3,5.539855593785816)-- (6.3,-5.539855593785816);
			\draw [line width=0.5pt] (6.8,6.102458520956943)-- (6.8,-6.102458520956942);

			\draw [line width=0.5pt] (7.3,6.655073252789935)-- (7.3,2.071231517720798);
			\draw [line width=0.5pt] (7.8,7.2)-- (7.8,3.44093010681705);
			\draw [line width=0.5pt] (8.3,7.738862965578342)-- (8.3,4.459820624195554);
			\draw [line width=0.5pt] (8.8,8.272847151978574)-- (8.8,5.332916650389355);
			\draw [line width=0.5pt] (9.3,8.802840450672726)-- (9.3,6.122907805936654);
			\draw [line width=0.5pt] (9.8,9.329523031752482)-- (9.8,6.8585712797929);
			\draw [line width=0.5pt] (10.3,9.853425800197616)-- (10.3,7.5557924799454375);
			\draw [line width=0.5pt] (10.8,10.37496987947435)-- (10.8,8.224354077980836);
			\draw [line width=0.5pt] (11.3,10.894494022211404)-- (11.3,8.870738413458037);
			\draw [line width=0.5pt] (11.8,11.412274094149685)-- (11.8,9.499473669630333);
			\draw [line width=0.5pt] (12.3,11.928537211242627)-- (12.3,10.113851887386923);
			\draw [line width=0.5pt] (18.3,18.05242365999646)-- (18.3,16.908281994336388);
			\draw [line width=0.5pt] (17.8,17.54536975956905)-- (17.8,16.36581803638303);
			\draw [line width=0.5pt] (17.3,17.03789893149974)-- (17.3,15.820556248122251);
			\draw [line width=0.5pt] (16.8,16.529972776747094)-- (16.8,15.27219696049);
			\draw [line width=0.5pt] (16.3,16.021547990128795)-- (16.3,14.720394016465727);
			\draw [line width=0.5pt] (15.8,15.512575543732254)-- (15.8,14.164744967700619);
			\draw [line width=0.5pt] (15.3,15.002999700059986)-- (15.3,13.604778572251737);
			\draw [line width=0.5pt] (14.8,14.492756811593853)-- (14.8,13.039938650162432);
			\draw [line width=0.5pt] (14.3,13.98177385026664)-- (14.3,12.469562943423478);
			\draw [line width=0.5pt] (13.8,13.46996659238619)-- (13.8,11.892854997854805);
			\draw [line width=0.5pt] (13.3,12.95723735986958)-- (13.3,11.308846094982458);
			\draw [line width=0.5pt] (12.8,12.443472184241825)-- (12.8,10.716342659695053);

			\draw [line width=0.5pt] (7.3,-6.655073252789935)-- (7.3,-2.071231517720798);
			\draw [line width=0.5pt] (7.8,-7.2)-- (7.8,-3.44093010681705);
			\draw [line width=0.5pt] (8.3,-7.738862965578342)-- (8.3,-4.459820624195554);
			\draw [line width=0.5pt] (8.8,-8.272847151978574)-- (8.8,-5.332916650389355);
			\draw [line width=0.5pt] (9.3,-8.802840450672726)-- (9.3,-6.122907805936654);
			\draw [line width=0.5pt] (9.8,-9.329523031752482)-- (9.8,-6.8585712797929);
			\draw [line width=0.5pt] (10.3,-9.853425800197616)-- (10.3,-7.5557924799454375);
			\draw [line width=0.5pt] (10.8,-10.37496987947435)-- (10.8,-8.224354077980836);
			\draw [line width=0.5pt] (11.3,-10.894494022211404)-- (11.3,-8.870738413458037);
			\draw [line width=0.5pt] (11.8,-11.412274094149685)-- (11.8,-9.499473669630333);
			\draw [line width=0.5pt] (12.3,-11.928537211242627)-- (12.3,-10.113851887386923);
			\draw [line width=0.5pt] (18.3,-18.05242365999646)-- (18.3,-16.908281994336388);
			\draw [line width=0.5pt] (17.8,-17.54536975956905)-- (17.8,-16.36581803638303);
			\draw [line width=0.5pt] (17.3,-17.03789893149974)-- (17.3,-15.820556248122251);
			\draw [line width=0.5pt] (16.8,-16.529972776747094)-- (16.8,-15.27219696049);
			\draw [line width=0.5pt] (16.3,-16.021547990128795)-- (16.3,-14.720394016465727);
			\draw [line width=0.5pt] (15.8,-15.512575543732254)-- (15.8,-14.164744967700619);
			\draw [line width=0.5pt] (15.3,-15.002999700059986)-- (15.3,-13.604778572251737);
			\draw [line width=0.5pt] (14.8,-14.492756811593853)-- (14.8,-13.039938650162432);
			\draw [line width=0.5pt] (14.3,-13.98177385026664)-- (14.3,-12.469562943423478);
			\draw [line width=0.5pt] (13.8,-13.46996659238619)-- (13.8,-11.892854997854805);
			\draw [line width=0.5pt] (13.3,-12.95723735986958)-- (13.3,-11.308846094982458);
			\draw [line width=0.5pt] (12.8,-12.443472184241825)-- (12.8,-10.716342659695053);

			\draw [line width=0.5pt] (-3.3,1.3747727084867516)-- (-3.3,-1.3747727084867516);
			\draw [line width=0.5pt] (-3.8,2.33238075793812)-- (-3.8,-2.3323807579381204);
			\draw [line width=0.5pt] (-4.3,3.080584360149872)-- (-4.3,-3.0805843601498726);
			\draw [line width=0.5pt] (-4.8,3.746998799039039)-- (-4.8,-3.7469987990390385);
			\draw [line width=0.5pt] (-5.3,4.369210454990696)-- (-5.3,-4.369210454990696);
			\draw [line width=0.5pt] (-5.8,4.963869458396343)-- (-5.8,-4.963869458396342);
			\draw [line width=0.5pt] (-6.3,5.539855593785816)-- (-6.3,-5.539855593785816);
			\draw [line width=0.5pt] (-6.8,6.102458520956943)-- (-6.8,-6.102458520956942);

			\draw [line width=0.5pt] (-7.3,6.655073252789935)-- (-7.3,2.071231517720798);
			\draw [line width=0.5pt] (-7.8,7.2)-- (-7.8,3.44093010681705);
			\draw [line width=0.5pt] (-8.3,7.738862965578342)-- (-8.3,4.459820624195554);
			\draw [line width=0.5pt] (-8.8,8.272847151978574)-- (-8.8,5.332916650389355);
			\draw [line width=0.5pt] (-9.3,8.802840450672726)-- (-9.3,6.122907805936654);
			\draw [line width=0.5pt] (-9.8,9.329523031752482)-- (-9.8,6.8585712797929);
			\draw [line width=0.5pt] (-10.3,9.853425800197616)-- (-10.3,7.5557924799454375);
			\draw [line width=0.5pt] (-10.8,10.37496987947435)-- (-10.8,8.224354077980836);
			\draw [line width=0.5pt] (-11.3,10.894494022211404)-- (-11.3,8.870738413458037);
			\draw [line width=0.5pt] (-11.8,11.412274094149685)-- (-11.8,9.499473669630333);
			\draw [line width=0.5pt] (-12.3,11.928537211242627)-- (-12.3,10.113851887386923);
			\draw [line width=0.5pt] (-18.3,18.05242365999646)-- (-18.3,16.908281994336388);
			\draw [line width=0.5pt] (-17.8,17.54536975956905)-- (-17.8,16.36581803638303);
			\draw [line width=0.5pt] (-17.3,17.03789893149974)-- (-17.3,15.820556248122251);
			\draw [line width=0.5pt] (-16.8,16.529972776747094)-- (-16.8,15.27219696049);
			\draw [line width=0.5pt] (-16.3,16.021547990128795)-- (-16.3,14.720394016465727);
			\draw [line width=0.5pt] (-15.8,15.512575543732254)-- (-15.8,14.164744967700619);
			\draw [line width=0.5pt] (-15.3,15.002999700059986)-- (-15.3,13.604778572251737);
			\draw [line width=0.5pt] (-14.8,14.492756811593853)-- (-14.8,13.039938650162432);
			\draw [line width=0.5pt] (-14.3,13.98177385026664)-- (-14.3,12.469562943423478);
			\draw [line width=0.5pt] (-13.8,13.46996659238619)-- (-13.8,11.892854997854805);
			\draw [line width=0.5pt] (-13.3,12.95723735986958)-- (-13.3,11.308846094982458);
			\draw [line width=0.5pt] (-12.8,12.443472184241825)-- (-12.8,10.716342659695053);

			\draw [line width=0.5pt] (-7.3,-6.655073252789935)-- (-7.3,-2.071231517720798);
			\draw [line width=0.5pt] (-7.8,-7.2)-- (-7.8,-3.44093010681705);
			\draw [line width=0.5pt] (-8.3,-7.738862965578342)-- (-8.3,-4.459820624195554);
			\draw [line width=0.5pt] (-8.8,-8.272847151978574)-- (-8.8,-5.332916650389355);
			\draw [line width=0.5pt] (9.3,-8.802840450672726)-- (9.3,-6.122907805936654);
			\draw [line width=0.5pt] (-9.8,-9.329523031752482)-- (-9.8,-6.8585712797929);
			\draw [line width=0.5pt] (-10.3,-9.853425800197616)-- (-10.3,-7.5557924799454375);
			\draw [line width=0.5pt] (-10.8,-10.37496987947435)-- (-10.8,-8.224354077980836);
			\draw [line width=0.5pt] (-11.3,-10.894494022211404)-- (-11.3,-8.870738413458037);
			\draw [line width=0.5pt] (-11.8,-11.412274094149685)-- (-11.8,-9.499473669630333);
			\draw [line width=0.5pt] (-12.3,-11.928537211242627)-- (-12.3,-10.113851887386923);
			\draw [line width=0.5pt] (-18.3,-18.05242365999646)-- (-18.3,-16.908281994336388);
			\draw [line width=0.5pt] (-17.8,-17.54536975956905)-- (-17.8,-16.36581803638303);
			\draw [line width=0.5pt] (17.3,-17.03789893149974)-- (17.3,-15.820556248122251);
			\draw [line width=0.5pt] (-16.8,-16.529972776747094)-- (-16.8,-15.27219696049);
			\draw [line width=0.5pt] (-16.3,-16.021547990128795)-- (-16.3,-14.720394016465727);
			\draw [line width=0.5pt] (-15.8,-15.512575543732254)-- (-15.8,-14.164744967700619);
			\draw [line width=0.5pt] (-15.3,-15.002999700059986)-- (-15.3,-13.604778572251737);
			\draw [line width=0.5pt] (-14.8,-14.492756811593853)-- (-14.8,-13.039938650162432);
			\draw [line width=0.5pt] (-14.3,-13.98177385026664)-- (-14.3,-12.469562943423478);
			\draw [line width=0.5pt] (-13.8,-13.46996659238619)-- (-13.8,-11.892854997854805);
			\draw [line width=0.5pt] (-13.3,-12.95723735986958)-- (-13.3,-11.308846094982458);
			\draw [line width=0.5pt] (-12.8,-12.443472184241825)-- (-12.8,-10.716342659695053);

			\pgftransformrotate{45}
			\pgftransformshift{\pgfpoint{11.9cm}{0.51cm}}
			\pgfnode{rectangle}{center}
			{\color{blue} \Huge{$\{g=0\}^+$}}{}{}
			
			\pgftransformrotate{-1}
			\pgftransformshift{\pgfpoint{7.9cm}{-2.51cm}}
			\pgfnode{rectangle}{center}
			{\color{red} \Huge{$\{h=0\}^+$}}{}{}

			\pgftransformrotate{-362}
			\pgftransformshift{\pgfpoint{-26.9cm}{0.3cm}}
			\pgfnode{rectangle}{center}
			{\color{blue} \Huge{$\{g=0\}^-$}}{}{}

			\pgftransformrotate{4}
			\pgftransformshift{\pgfpoint{-3.5cm}{3.8cm}}
			\pgfnode{rectangle}{center}
			{\color{red} \Huge{$\{h=0\}^-$}}{}{}

			\pgftransformrotate{315}
			\pgftransformshift{\pgfpoint{7cm}{15cm}}
			\pgfnode{rectangle}{center}
			{\color{black} \Huge{$\{f_1=23.04\}$}}{}{}

			\pgftransformrotate{0}
			\pgftransformshift{\pgfpoint{19cm}{-5cm}}
			\pgfnode{rectangle}{center}
			{\color{black} \Huge{$\{f_2=25.5\}$}}{}{}

			\end{tikzpicture}
		\end{center}
		\caption{Graphic representation for Example \ref{24e}}
		\label{xxsdaq}
	\end{figure}
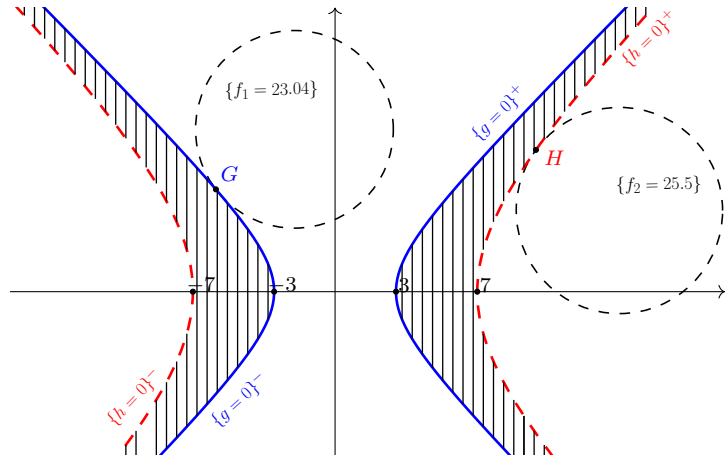

	\begin{example}\label{25e}	The following two problems belong to [Non-Alter]. Both Assumptions 1 and 2 are satisfied and the level sets are not arranged alternatively.  They can be solved in polynomial time for any given quadratic objective function.
		
		- $g(x,y)=x^2+3y^2-16,~ h(x,y)=-2x^2-y^2+4.$ (the left picture in Figure 4)
		
		- $g(x,y)=-x^2+3y^2+4,~ h(x,y)=y^2-x-2$ (the right picture in Figure 4)
		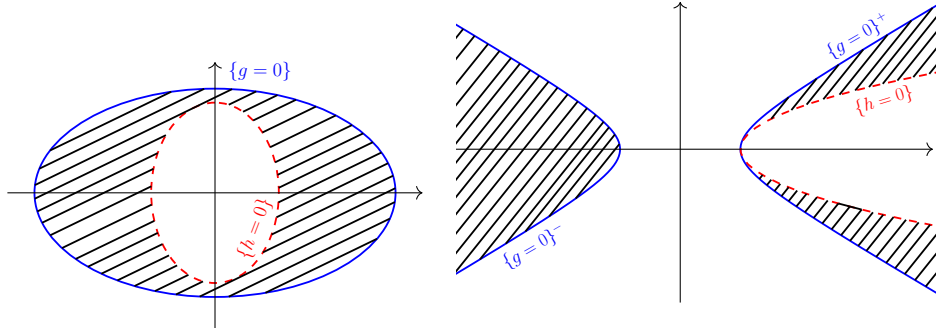
\begin{figure}	\begin{center} 	
				\begin{tikzpicture}[scale=0.6]
				\clip(-5.740246485385321,-2.996805788698608) rectangle (5.129711328425287,2.927427579493658);
				\draw [rotate around={0.:(0.,0.)},line width=0.7pt,color=blue] (0.,0.) ellipse (4.cm and 2.309401076758503cm);
				\draw [rotate around={90.:(0.,0.)},line width=0.7pt,color=red,dash pattern=on 3pt off 3pt] (0.,0.) ellipse (2.cm and 1.4142135623730951cm);
				\draw [line width=0.7pt] (-3.6544224222135027,0.9390059921752711)-- (-0.9702215324074972,2.2404367266266676);
				\draw [line width=0.7pt] (-3.891095237139262,0.5352189133155874)-- (-0.24045877558539283,2.3052244704326137);
				\draw [line width=0.7pt] (-3.9935840255915602,-0.13074992739311075)-- (-1.0877370843433387,1.278145559272693);
				\draw [line width=0.7pt] (-3.988925058993888,0.17173319784926872)-- (0.37822140914523844,2.2990540493620553);
				\draw [line width=0.7pt] (0.3043997393245416,1.9531209889298466)-- (0.913385160041071,2.248386647459073);
				\draw [line width=0.7pt] (1.3834395393682728,2.1668790952975434)-- (0.6207039911964697,1.7970679204263664);
				\draw [line width=0.7pt] (-1.2698014770130888,0.880459208567187)-- (-3.935787326892748,-0.41213999137446544);
				\draw [line width=0.7pt] (-3.8311486913271104,-0.6639025292271363)-- (-1.3631240231333572,0.532715491715289);
				\draw [line width=0.7pt] (-3.6519958857824153,-0.9421475557166031)-- (-1.4105107699886743,0.14463310648642236);
				\draw [line width=0.7pt] (-1.402103672709185,-0.26117155655015206)-- (-3.390582371769111,-1.2252824409428429);
				\draw [line width=0.7pt] (-3.088781904515843,-1.467358891152795)-- (-1.3441981987376999,-0.6215001247149078);
				\draw [line width=0.7pt] (-1.2671411329628715,-0.8880916046822764)-- (-2.817024338877535,-1.639550128762113);
				\draw [line width=0.7pt] (-2.518265386807012,-1.79428160569532)-- (-1.1632912559576076,-1.1373244513440939);
				\draw [line width=0.7pt] (-1.0183656148161586,-1.387754642982763)-- (-2.162884175561641,-1.9426727330411784);
				\draw [line width=0.7pt] (-0.8391177493266936,-1.6098952778146187)-- (-1.7826666128393853,-2.067373514669257);
				\draw [line width=0.7pt] (-0.5882339862694995,-1.818780238180245)-- (-1.329388835233008,-2.1781280437383095);
				\draw [line width=0.7pt] (-0.2659860639079383,-1.9643072131450123)-- (-0.8650291677380301,-2.2547523543959658);
				\draw [line width=0.7pt] (-0.33896279704765925,-2.3010942485563706)-- (3.9842076668168454,-0.2050115994099453);
				\draw [line width=0.7pt] (0.3049386392983674,-2.3026804979025974)-- (3.875609323716022,-0.5714462266697959);
				\draw [line width=0.7pt] (3.5898131862570004,-1.0186987267711294)-- (1.1173438522468702,-2.217471737200283);
				\draw [line width=0.7pt] (0.8452657408426326,1.60344992273395)-- (1.7948529577315964,2.0638558460740537);
				\draw [line width=0.7pt] (1.051055521556521,1.3381197932962134)-- (2.238564424257597,1.9138816855149168);
				\draw [line width=0.7pt] (1.2204767554083873,1.010382590415943)-- (2.6763983383116927,1.7162839639448177);
				\draw [line width=0.7pt] (1.3310817870471658,0.6756053229383598)-- (3.0372062242981386,1.5028171713024672);
				\draw [line width=0.7pt] (1.3858956417827129,0.3982292558004304)-- (3.2839847176618555,1.3185148683478929);
				\draw [line width=0.7pt] (1.4120215343343816,0.11131205303990033)-- (3.4963098093741665,1.1218760651804014);
				\draw [line width=0.7pt] (1.4062315590248369,-0.21219237687803838)-- (3.6880307154710503,0.8941344868534585);
				\draw [line width=0.7pt] (1.3610656234914305,-0.5431397031149228)-- (3.8350432295065437,0.6563645907105862);
				\draw [line width=0.7pt] (1.253581615911387,-0.9257787340925427)-- (3.945634603485601,0.3794590780646514);
				\draw [line width=0.7pt] (1.0622570385834798,-1.3203105573915979)-- (3.996084987383052,0.10215147839001237);
				\draw [->] (-4.6,0)--(4.6,0);
				\draw [->] (0,-5.1)--(0,2.9);
				
				
				
				\pgftransformrotate{0}
				\pgftransformshift{\pgfpoint{1cm}{2.7cm}}
				\pgfnode{rectangle}{center}
				{\color{blue} \normalsize{$\{g=0\}$}}{}{}
				
				\pgftransformrotate{69}
				\pgftransformshift{\pgfpoint{-3.3cm}{-1.146cm}}
				\pgfnode{rectangle}{center}
				{\color{red} \normalsize{$\{h=0\}$}}{}{}
				
				\end{tikzpicture}
				\begin{tikzpicture}[scale=0.4]
				\clip(-7.419809524079226,-5.942145651756148) rectangle (8.470141793728978,5.47980666971277);
				\draw [samples=50,domain=-0.99:0.99,rotate around={0.:(0.,0.)},xshift=0.cm,yshift=0.cm,line width=0.7pt,color=blue] plot ({2.*(1+(\x)^2)/(1-(\x)^2)},{1.1547005383792515*2*(\x)/(1-(\x)^2)});
				\draw [samples=50,domain=-0.99:0.99,rotate around={0.:(0.,0.)},xshift=0.cm,yshift=0.cm,line width=0.7pt,color=blue] plot ({2.*(-1-(\x)^2)/(1-(\x)^2)},{1.1547005383792515*(-2)*(\x)/(1-(\x)^2)});
				\draw [samples=50,rotate around={-90.:(2.,0.)},xshift=2.cm,yshift=0.cm,color=red,line width=0.7pt,dash pattern=on 3pt off 3pt, domain=-5.0:5.0)] plot (\x,{(\x)^2/2/0.5});
				\draw [line width=0.7pt] (-2.1723428404324947,0.48958261692259136)-- (-4.35710736601543,-2.234903771156297);
				\draw [line width=0.7pt] (-2.5476026733303994,0.9110944665913627)-- (-5.797587066606141,-3.1417625517575747);
				\draw [line width=0.7pt] (-3.098597279000436,1.3664192984376977)-- (-7.578716601517537,-4.220463849946572);
				\draw [line width=0.7pt] (-3.7416092313872493,1.825708961874107)-- (-9.510166876893631,-5.3679069175475576);
				\draw [line width=0.7pt] (-4.375685919740274,2.246970350510014)-- (-11.350176920102065,-6.450491353832865);
				\draw [line width=0.7pt] (-4.982982919000477,2.635028068860654)-- (-12.017868115860832,-6.13774754888229);
				\draw [line width=0.7pt] (3.599892030121063,1.7281225485604454)-- (3.026610133199264,1.0132177126359683);
				\draw [line width=0.7pt] (5.267237412219601,2.8132892940877254)-- (4.200919795452928,1.4835497280013663);
				\draw [line width=0.7pt] (3.025712565051508,-1.0127746862217222)-- (2.8790460141420517,-1.1956735830411291);
				\draw [line width=0.7pt] (3.601522383582312,-1.265512695938809)-- (3.364018487488721,-1.5616892973279395);
				\draw [line width=0.7pt] (4.19816862050061,-1.4826222109831653)-- (3.8590823327923007,-1.905476008006894);
				\draw [line width=0.7pt] (4.902244138577569,-1.7035974109447245)-- (4.435456340076356,-2.2856999616983518);
				\draw [line width=0.7pt] (5.63747044708653,-1.9072153646315173)-- (5.030061740068232,-2.664677660862458);
				\draw [line width=0.7pt] (6.405431504887387,-2.098911981214883)-- (5.644714068414617,-3.047556229405122);
				\draw [line width=0.7pt] (7.089646564335225,-2.2560245043738383)-- (6.187740623849113,-3.380736311714388);
				\draw [line width=0.7pt] (7.924956683583716,-2.434123391199328)-- (6.845737093638954,-3.7799522199460096);
				\draw [line width=0.7pt] (8.714116329425842,-2.5911611932540675)-- (7.4631043759070295,-4.1512217850364115);
				\draw [line width=0.7pt] (6.890812588948338,3.807155462021871)-- (5.293042739466973,1.8146742791660913);
				\draw [line width=0.7pt] (8.38758262044066,4.702890678958996)-- (6.2741935722251805,2.067412288883178);
				\draw [line width=0.7pt] (9.854573037431393,5.571134227727474)-- (7.219039872620316,2.2845218039275346);
				\draw [line width=0.7pt] (11.521622456948545,6.551024959003515)-- (8.277515236497225,2.5054970038890936);
				\draw [line width=0.7pt] (13.214307269778464,7.541395021703432)-- (9.339303853361395,2.7091149575758857);
				\draw [line width=0.7pt] (14.132918325956636,7.537564815956076)-- (10.41470778877628,2.9008115741592526);
				\draw [line width=0.7pt] (15.079837391450504,7.7080544120118795)-- (11.350899784959378,3.057924097318208);
				\draw [line width=0.7pt] (16.04110914897193,7.687034590173117)-- (12.471844753906277,3.2360229841436965);
				\draw [line width=0.7pt] (17.063651386632305,7.821033962686987)-- (13.512861498837552,3.3930607861984363);
				\draw [line width=0.7pt] (-5.622697431546987,3.033959261796653)-- (-13.28608658942122,-6.52258525487561);
				\draw [line width=0.7pt] (-14.294620567980765,-6.360913198305853)-- (-6.385983148842973,3.5014654445106466);
				\draw [line width=0.7pt] (-7.213460358241191,4.001416846564774)-- (-14.915029426207427,-5.602739503843204);
				\draw [line width=0.7pt] (16.74904322297811,6.1466147567628955)-- (14.674995013362228,3.5601959234517175);
				\draw [line width=0.7pt] (9.608198646890305,-2.7582963305073487)-- (8.158250885719857,-4.566437616433948);
				\draw [line width=0.7pt] (2.736088346403543,-0.857955911689839)-- (2.6319663662480126,-0.9878000730682241);
				\draw [line width=0.7pt] (3.8862899702898415,-1.3734227209019958)-- (3.601153800624853,-1.7289986211441215);
				\draw [line width=0.7pt] (3.325702861803686,-1.1513917065029111)-- (3.1326879499164018,-1.3920888132049893);
				\draw [line width=0.7pt] (4.545288865515151,-1.5953961468911568)-- (4.144193401431107,-2.095577800072105);
				\draw [line width=0.7pt] (4.721736215284572,-2.4694663989567918)-- (5.2551993548895535,-1.8042171030365366);
				\draw [line width=0.7pt] (5.2551993548895535,-1.8042171030365366)-- (6.029957906735843,-2.0074755058868945);
				\draw [line width=0.7pt] (6.029957906735843,-2.0074755058868945)-- (5.34493787526335,-2.86172214411909);
				\draw [line width=0.7pt] (6.746826135503297,-2.1787212156453837)-- (5.9161571779365705,-3.214597732949407);
				\draw [line width=0.7pt] (7.504133958862277,-2.346089077350278)-- (6.514871973119891,-3.5797372727684476);
				\draw [line width=0.7pt] (8.305137989652977,-2.511003383042917)-- (7.143635937446912,-3.959441644423407);
				\draw [line width=0.7pt] (9.151031030700473,-2.6741411762845417)-- (7.8033358956156595,-4.354769458107735);
				\draw [line width=0.7pt] (18.106227522878495,8.493337359646421)-- (14.08285942934152,3.476040769228911);
				\draw [line width=0.7pt] (17.35887494988688,8.779358951461962)-- (12.97532612830542,3.312902975987287);
				\draw [line width=0.7pt] (16.288136866660604,8.607904396503132)-- (11.909832668303459,3.1479886702946467);
				\draw [line width=0.7pt] (15.346105245436188,8.54491449335722)-- (10.88410040459823,2.9806208085897525);
				\draw [line width=0.7pt] (14.10323591553502,8.060216359569944)-- (9.892588445933171,2.809375098831263);
				\draw [line width=0.7pt] (12.33924587372621,7.029864596426433)-- (8.791844233070433,2.6061166959809063);
				\draw [line width=0.7pt] (10.683322621392707,6.0589708210982804)-- (7.747026864233564,2.3972957398355264);
				\draw [line width=0.7pt] (9.095650597102498,5.122852713888069)-- (6.732027169117903,2.1753223138463644);
				\draw [line width=0.7pt] (7.6838097426913095,4.283337178800267)-- (5.815346900496445,1.9532912994472802);
				\draw [line width=0.7pt] (6.084567996972266,3.3177285256522064)-- (4.755120296264482,1.659855504634208);
				\draw [line width=0.7pt] (4.53526067936403,2.350077263253663)-- (3.6946316919214377,1.3017802010790598);
				\draw [line width=0.7pt] (-6.8339032128968356,3.7728076690638797)-- (-14.960415306811242,-6.361269172247317);
				\draw [line width=0.7pt] (-6.005609111900549,3.2694209887451993)-- (-13.843140886874792,-6.504286148574625);
				\draw [line width=0.7pt] (-5.324657023414291,2.849091809506738)-- (-13.265603193783052,-7.05274858046961);
				\draw [line width=0.7pt] (-13.265603193783052,-7.05274858046961)-- (-12.715641140330273,-6.754197751452386);
				\draw [line width=0.7pt] (-4.685814563669582,2.4465525490029125)-- (-11.490011421206933,-6.062816884254604);
				\draw [line width=0.7pt] (-11.490011421206933,-6.062816884254604)-- (-11.819988653278601,-6.361367713271828);
				\draw [line width=0.7pt] (-4.076311313718056,2.0507164541809795)-- (-10.486822166272994,-5.943440057380759);
				\draw [line width=0.7pt] (-3.405064409558854,1.5910650136773563)-- (-8.510991873104107,-4.776225241200564);
				\draw [line width=0.7pt] (-6.658209858176487,-3.666595083339971)-- (-2.8064567075270195,1.136690994542697);
				\draw [line width=0.7pt] (-2.366167701866162,0.7300113225958881)-- (-5.147495289843418,-2.738412420786916);
				\draw [line width=0.7pt] (-2.021948049960682,0.17153611547796016)-- (-3.475401503443869,-1.6409769051939667);
				
				\draw [->] (-7.6,0)--(8.4,0);
				\draw [->] (0,-5.1)--(0,4.9);
				
				\pgftransformrotate{30}
				\pgftransformshift{\pgfpoint{7cm}{0.3cm}}
				\pgfnode{rectangle}{center}
				{\color{blue} \Large{$\{g=0\}^+$}}{}{}
				
				\pgftransformrotate{-13}
				\pgftransformshift{\pgfpoint{0.2cm}{-2.3cm}}
				\pgfnode{rectangle}{center}
				{\color{red} \Large{$\{h=0\}$}}{}{}
				
				\pgftransformrotate{13}
				\pgftransformshift{\pgfpoint{-12.5cm}{1.6cm}}
				\pgfnode{rectangle}{center}
				{\color{blue} \Large{$\{g=0\}^-$}}{}{}

				\end{tikzpicture}
				\caption{Graphic representation for Example \ref{25e}}
				\label{25e-fig}
			\end{center}
		\end{figure}
		
	\end{example}

	Finally in this section, we show that ${\rm [YZ]} \subset \hbox{[Non-Alter]}.$
	
	For (P) to be in [YZ], it is required that (i) the SDP relaxation (SP) and
	its conic dual (SD) both have a Slater point; (ii) $\{g=0\}\cap\{h=0\}=\emptyset;$ and (iii)
	$\{h\ge0\} \subset \{g\le0\}.$ We point out that
	(iii) implies parenthetically that
	\begin{equation}\label{YZ3}
	\{h\ge0\} \subset \{g\le0\}~\Rightarrow~\{g>0\} \subset \{h<0\}~\Rightarrow~
	\{g\ge0\} \subset \{h\le0\}.
	\end{equation}
	By requiring (iii), it assumes both $\{h\ge0\} \subset \{g\le0\}$ and
	$\{g\ge0\} \subset \{h\le0\},$ which shows that Assumption 2:
	$\{h=0\} \subset \{g\le0\}$ and $\{g=0\} \subset \{h\le0\}$ is a stronger version of (iii).
	Since [Non-Alter] does not assume the Slater condition; does not restrict from a common root of $g$ and $h,$ it follows that ${\rm [YZ]} \subset \hbox{[Non-Alter]}.$
	
	The SDP relaxation (SP) and its conic dual (SD) in [YZ] have the following formulations.
	\[\begin{array} { c l }  ( \rm SP )\,\,\,\,\,  \min & {M(f)\bullet X} \\ { \text { subject to  } } & { M(g)\bullet X \leq 0, } \\ { } & { M(h)\bullet X \leq 0,} \\ { } & { X_{00}=1, X\succeq 0};\end{array}\]
	\[\begin{array} { c l }  ( \rm SD )\,\,\,\,\,  \max & {\gamma} \\ { \text { subject to  } } & { Z=M(f)-\gamma \begin{bmatrix}1& \bar{0}^T \\\bar{0} & [0] \end{bmatrix}+\alpha M(g)+\beta M(h), } \\ { } & { Z\succeq 0, \alpha\geq 0, \beta\geq 0},\end{array}\]
	where $\bar{0}$ is zero vector in $\mathbb{ R }^n,$~$[0]$ is the zero matrix in $\mathbb{ R }^{n\times n}$, and
	\begin{equation}\label{M-notation}
	M(f)=\begin{bmatrix} a_0& a^T \\a & A \end{bmatrix};\,M(g)=\begin{bmatrix} b_0& b^T \\b & B \end{bmatrix};\,M(h)=\begin{bmatrix} c_0& c^T \\c & C \end{bmatrix}
	\end{equation}
	
	%
	For the rest of the paper, we prove that (SP) is tight for a larger subclass [Non-Alter] of (P). In particular, the deletion of the dual Slater condition
	\begin{equation}\label{Dual-Slater}
	(\exists\tilde{\alpha}\ge0,\tilde{\beta}\ge0)~M(f)-\gamma \begin{bmatrix}1& \bar{0}^T \\\bar{0} & [0] \end{bmatrix}+\tilde{\alpha} M(g)+\tilde{\beta} M(h)\succ0
	\end{equation}
	creates a notable technical gap to conquer. Without \eqref{Dual-Slater}, we do not have a ``full-dimensional'' dual to work with, and there is no easy way to returning back to a full-dimensional case, either. As such, strong duality for the nonconvex (P) becomes very tricky to establish. On the other hand, checking \eqref{Dual-Slater} is another subtle issue. Although the strict feasibility of
	\eqref{Dual-Slater} with two variables $\tilde{\alpha},~\tilde{\beta}$ can be determined in polynomial time by computational complexity theory \cite[1997]{PK97}, practical checkable algorithms seem to be unavailable so far.
	
	The analysis in the following section is fundamental in its own right.
	It concerns with the unsolvability of a system of two quadratic inequalities.

	\section{Unsolvability of \eqref{Unsolvable-system}}
	
	This section gives a detail account of the conditions under which \eqref{Unsolvable-system} has no solution, which will be shown to be a hinge for strong duality of (P) in next section.
	
	The unsolvability of a single quadratic inequality, for example $f(x)<0$, has a well-known result:
	
	$f(x)<0$ has no solution if and only if $\begin{pmatrix}A & a\\a^T & a_0 \end{pmatrix}\succeq 0.$
	
	\noindent For two quadratic inequalities, the situation is very different. The type of problem has not been studied much in literature. We first note that,
	\begin{eqnarray}
	&&\eqref{Unsolvable-system} \hbox{ is unsolvable }\label{tty}\\
	&\Longleftrightarrow& \Big(\{h\ge0\} \subset \{g\le0\}\Big)\wedge \Big(\{g\ge0\} \subset \{h\le0\})\Big) \nonumber\\
	&\Longrightarrow& \hbox{ Assumption 2 } ~~(\hbox{ see }\eqref{YZ3}) \label{tty-1}
	\end{eqnarray}
	The main result in this section, Theorem \ref{lm31aa} below, concerns the reverse from \eqref{tty-1} to \eqref{tty} under additional hypotheses.
	
	\begin{theorem}\label{lm31aa}
		Under Assumptions $1,~2,~3,~4$ and $5$, neither of the following two quadratic systems
		has a solution $($the unsolvability of \eqref{Unsolvable-system}$).$	
		$$  \begin{cases}
		\begin{array}{l}
		g(x)> 0,\\
		h(x)\geq 0;
		\end{array}\end{cases} \hbox{ and }\ \ \  \begin{cases}
		\begin{array}{l}
		g(x)\ge 0,\\
		h(x) >0.
		\end{array}\end{cases}$$
		In other words, $\big(\{g> 0\}\cap\{h\geq 0\}\big)
		\cup\big(\{g \geq 0\}\cap\{h> 0\}\big)=\emptyset.$	
	\end{theorem}
	
	\noindent {\bf Assumption 3:} $f(x), g(x), h(x)$ are non-constant functions and $$\mathcal{D}=\{x|g(x)\leq 0, h(x)\leq 0 \}\ne \emptyset; ~ \mathcal{D}\ne \{g\leq 0\};~ \mathcal{D}\ne \{h\leq 0\}.$$
	
	\noindent {\bf Assumption 4:} All the functions $g, h, -g, -h$ satisfy Slater's condition.
	
	\noindent {\bf Assumption 5:} Suppose, after coordinate transformation,
	$g,~h$ degenerates to single variable functions (in $x_1\in\mathbb{R}$) of the form $-x_1^2+1$ and ${\bar{c}_1}x_1+{\bar{c}_0},~\bar{c}_1>0$ one for each. Then, $\bar{c}_0\not=\pm \bar{c}_1.$

	Assumption 3 is a prior condition to avoid triviality.
	Assumption 4 is a regularity condition for the constraint functions $g,~h$ not to have a degenerate gradient at the boundary. Assumption 5 excludes one peculiar example: $g(x)=-x^2+1$ and $h(x)=x-1.$ It implies that $\{g=0\}\subset\{h\le0\},$ $\{h=0\}\subset \{g\le0\}$ but $\{h=0\}$ separates $\{g<0\}.$ Such a pathological case, satisfying Assumption 2 as well as the two-sided Slater condition Assumption 4, fails surprisingly the $\mathcal{S}$-lemma with equality \cite[2016]{Xia-Wang-Sheu16}\cite[2018]{Quang-Sheu18}.
	
	The main tools to prove the un-solvability result, Theorem \ref{lm31aa}, are ``the separation theorem'' and ``the $\mathcal{S}$-lemma with equality'' of Quang and Sheu \cite[2018]{Quang-Sheu18}, which we take, respectively, as Lemma \ref{th23sss} and Lemma \ref{th23sss-1} below. While the well-known separation theorem asserts the existence of a hyperplane to separate two disjoint convex sets, our
	theorem in \cite[2018]{Quang-Sheu18} may be thought of as a nonlinear version of the classical separation theorem. It describes when and how a hyper-surface $\{h=0\}$ can separate two convex connected branches of $\{g<0\}.$ Moreover, when the separation happens, the $\mathcal{S}$-lemma with equality fails.

	\begin{lemma}{\rm(Quang and Sheu \cite[2018]{Quang-Sheu18})}\label{th23sss}
		The hypersurface $\{h = 0\}$ separates $\{g < 0\}$ if and only if there exists a basis of $\mathbb{R}^n$ such that
		\begin{itemize}
			\item[$\rm(i)$] $g(x)$ is of the form $-x_1^2+\delta(x_2^2+\cdots
			+x_m^2)+\theta,~ \delta, \theta \in\{0, 1\};$
			\item[$\rm(ii)$] With the same basis and $\delta$ as in $\rm(i)$,~ $h(x)$ has the form $$\bar{c}_1x_1+\delta(\bar{c}_2x_2+\cdots
			+\bar{c}_mx_m)+\bar{c}_0,~ \bar{c}_1\ne 0;$$
			\item[$\rm(iii)$]
			$g|_{\{h = 0\}}(x)=-(\delta\dfrac{\bar{c}_2}{\bar{c}_1}x_2+\cdots+\delta\dfrac{\bar{c}_m}{\bar{c}_1}x_m+
			\dfrac{\bar{c}_0}{\bar{c}_1})^2+\delta(x_2^2+\cdots
			+x_m^2)+\theta\geq 0,~\forall (x_2, \cdots,
			x_n)^T\in\mathbb{R}^{n-1}.$
		\end{itemize}
	\end{lemma}

	\begin{lemma}{\rm(Quang and Sheu \cite[2018]{Quang-Sheu18})}\label{th23sss-1}
		If $h(x)$ takes both positive and negative values, then $\{h=0\}$ fails to separate $\{g<0\}$ if and only if
		$$\big(h(x)= 0~\Longrightarrow~ g(x)\ge0\big)\Longrightarrow (\exists \lambda\in\mathbb R)(\forall x\in\mathbb R^n)~g(x) + \lambda h(x)\ge0.$$
	\end{lemma}
	
	Before proceeding to prove Theorem \ref{lm31aa}, two technical lemmas, Lemma \ref{lemma3-3} and Lemma \ref{l24dd} will be shown in sequence to help.
	
	\begin{lemma}\label{lemma3-3}
		Under Assumptions 4 and 5, if either $\{g=0\}\subset \{h\leq 0\}$ or $\{g=0\}\subset \{h\geq 0\}$ holds true, $\{h=0\}$ separates neither $\{g<0\}$ nor $\{g>0\}.$ Similarly, if either $\{h=0\}\subset \{g\leq 0\}$ or $\{h=0\}\subset \{g\geq 0\}$ holds true, $\{g=0\}$ separates neither $\{h<0\}$ nor $\{h>0\}.$
	\end{lemma}
	\proof
	Suppose contrarily that $\{h=0\}$ separates $\{g<0\}$. We will show that both $\{g=0\}\subset \{h\leq 0\}$ and $\{g=0\}\subset \{h\geq 0\}$ are violated.
	
	From Lemma \ref{th23sss}, $g(x)$ is of the form
	$$-x_1^2+\delta(x_2^2+\cdots+x_m^2)+\theta,~ \delta, \theta \in\{0, 1\}.$$
	We first consider the case when rank$(B)=1,$ and then when rank$(B)>1.$
	
	When rank$(B)=1,$ there must be $\delta=0$ so that
	$g(x)=-x_1^2+\theta,~\theta \in\{0, 1\};~~h(x)=\bar{c}_1x_1+\bar{c}_0,~{\bar{c}_1}>0.$
	If $\theta=0,$ $g(x)\le0, ~\forall x,$ which fails the Slater condition Assumption 4. So, $g(x)=-x_1^2+1$ and
	$$\{g=0\}=\{-1,1\};~\{g<0\}=(-\infty,-1)\cup(1,\infty); ~\{g>0\}=(-1,1);$$
	$$\{h=0\}=-\frac{\bar{c}_0}{\bar{c}_1};~ \{h<0\}=(-\infty,-\frac{\bar{c}_0}{\bar{c}_1});~~\{h>0\}=(-\frac{\bar{c}_0}{\bar{c}_1},\infty).$$
	Under Assumption 5,  $-\frac{\bar{c}_0}{\bar{c}_1}\not=\pm 1.$

	- If $-\frac{\bar{c}_0}{\bar{c}_1}\in(-1,1),$ which contradicts to both $\{g=0\}\subset\{h\le0\}$ and $\{g=0\}\subset\{h\ge0\}.$
	
	- If $-\frac{\bar{c}_0}{\bar{c}_1}\in (-\infty,-1)\cup(1,\infty),$ then $\{h=0\}$  cannot separate $\{g<0\}.$

	In other words, rank$(B)=1$ is impossible under current circumstance.
	
	Now, suppose rank$(B)>1.$ Consider, for any $t\in\mathbb{R},$
	$$x(t)=(\sqrt{t^2+\theta}, -t,0, \cdots, 0)\in \{g=0\};~y(t)=(-\sqrt{t^2+\theta}, t, 0, \cdots, 0)\in \{g=0\}.$$
	Then, for any properly chosen $t_0$ such that ${\bar{c}_1}\sqrt{t_0^2+\theta}-\bar{c}_2t_0\ne 0,$ there is
	$$h(x(t_0))h(y(t_0))=-(\bar{c}_1\sqrt{t_0^2+\theta}-\bar{c}_2t_0)^2<0,$$
	which fails both $\{g=0\}\subset \{h\leq 0\}$ and $\{g=0\}\subset \{h\geq 0\}$.
	
	In summary, we have shown that $\{h=0\}$ cannot separate $\{g<0\}$ if either $\{g=0\}\subset \{h\leq 0\}$ or $\{g=0\}\subset \{h\geq 0\}$ holds true. The other parts of the lemma can be similarly argued.
	\hfill$\Box$

	In the next lemma, we use a variable transformation to make $g(x)$ adopt
	one of the following five canonical forms (see Quang and Sheu \cite[2018]{Quang-Sheu18}):
	\begin{eqnarray}
	&& -x_1^2-\cdots - x_k^2+\delta (x_{k+1}^2+\cdots +
	x_m^2)+\theta; \label{form:1}\\
	&&-x_1^2-\cdots - x_k^2+\delta(x_{k+1}^2+\cdots + x_m^2)-1;\label{form:2}\\
	&&-x_1^2-\cdots - x_k^2+\delta (x_{k+1}^2+\cdots+ x_m^2)+x_{m+1};\label{form:3}\\
	&&\hskip8pt x_1^2+\cdots + x_m^2+\eta x_{m+1}+c';\label{form:4}\\
	&& \hskip8pt \eta x_1+c', \label{form:5}
	\end{eqnarray}where $\delta, \eta, \theta\in \{0, 1\}.$
	
	\begin{lemma}\label{l24dd} Assume that $\{g=0\}=\{h=0\}$ and Assumptions 3 and 4 hold true.  Let $g(\bar{x})>0, h(\bar{x})>0$ for some $\bar{x}\in\mathbb{R}^n.$ Then, there exists a largest open ball
		$\mathcal{B}(\bar{x}, \bar{r})=\{x: \|x-\bar{x}\|<\bar{r}\}$ centered at $\bar{x}$ with the radius $\bar{r}$ such that
		\begin{equation}\label{xcxc}
		\mathcal{B}(\bar{x}, \bar{r})\subset \{g>0\}\cap \{h>0\} \text{ and } \overline{\mathcal{B}(\bar{x}, \bar{r})}\cap \{g=0\}\ne \emptyset. 
		\end{equation}	
		Moreover, $\forall x^*\in \overline{\mathcal{B}(\bar{x}, \bar{r})}\cap\{g=0\} ,$ $$\nabla g(x^*)=t^*\overrightarrow{x^*\bar{x}}\ne 0 \hbox{ and }~\nabla h(x^*)=s^*\overrightarrow{x^*\bar{x}}\ne 0,~~t^*>0,~s^*>0.$$
	\end{lemma}
	
	\proof Since $g(\bar{x})>0$ and $h(\bar{x})>0$, it follows that  $\{g>0\}\cap \{h>0\}$ is a non-empty open set. Hence the set $\mathcal{R}=\{r\in \mathbb{R}_+: \mathcal{B}(\bar{x}, {r})\subset \{g>0\}\cap \{h>0\}\}$ is non-empty. Let $\bar{r}=\sup\{r: r\in\mathcal{R}\}.$ Then, $\bar{r}<\infty$. Otherwise, $\bar{r}=\infty$ would imply that $\mathcal{B}(\bar{x}, {r})=\mathbb{ R }^n$ and $\mathcal{D}= \emptyset,$ which is a contradiction to Assumption 3. We therefore have shown that there exists a largest open ball
	$\mathcal{B}(\bar{x}, \bar{r})=\{x: \|x-\bar{x}\|<\bar{r}\}$ centered at $\bar{x}$ with the radius $\bar{r}$ such that
	$\mathcal{B}(\bar{x}, \bar{r})\subset \{g>0\}\cap \{h>0\}$.

It is clear that, if $x\in \left(\{g>0\}\cap \{h>0\}\right)^c,$ $\|x-\bar x\|\ge\bar r.$ Therefore, 
\begin{equation}\label{27kl-1}
\inf_{x\in\left(\{g>0\}\cap \{h>0\}\right)^c}\|x-\bar x\|\ge\bar r.
\end{equation}
On the other hand, since $\bar{r}$ is the largest radius such that $\mathcal{B}(\bar{x}, \bar{r})\subset \{g>0\}\cap \{h>0\}$, $\forall\varepsilon>0,$ there is some $x\in \left(\{g>0\}\cap \{h>0\}\right)^c$ such that $\|x-\bar x\|<\bar r+\varepsilon.$ Then,
\begin{equation}\label{27kl-2}
\inf_{x\in\left(\{g>0\}\cap \{h>0\}\right)^c}\|x-\bar x\|\le\bar r.
\end{equation}
Combining \eqref{27kl-1} and \eqref{27kl-2}, one has
\begin{equation}\label{27uuu2}
	\bar{r}^2 = \inf_{x\in\left(\{g>0\}\cap \{h>0\}\right)^c}\|x-\bar x\|^2= \min\left\{\inf_{\{g\le0\}}\|x-\bar x\|^2, \inf_{\{h\le0\}}\|x-\bar x\|^2\right\}.
	\end{equation}
Suppose, without loss of generality, in \eqref{27uuu2} there is $\bar{r}^2 = \inf_{\{g\le0\}}\|x-\bar x\|^2.$
Since the distance from a point to a closed set is always attainable, there is an $\hat{x}\in \{g\le 0\}$ such that 
\begin{equation*}
		\|\hat{x}-\bar{x}\|^2=\bar{r}^2,~i.e.,~ \hat{x}\in\overline{\mathcal{B}(\bar{x}, \bar{r})}.
\end{equation*}
Moreover, $\hat{x}$ cannot be an interior point of $\{g\le0\},$ for otherwise the distance from $\bar{x}$ to $\{g\le0\}$ would have been shorter than $\bar{r}$ (since $\mathcal{B}(\bar{x}, \bar{r})\subset \{g>0\}$).
This implies that $g(\bar{x})=0$ and thus 
$$\hat{x}\in\overline{\mathcal{B}(\bar{x}, \bar{r})}\cap\{g=0\}\ne \emptyset. $$


Now, $\forall x^*\in\overline{\mathcal{B}(\bar{x}, \bar{r})}\cap \{g=0\}$, since $x^*\not\in \{g>0\}\cap\{h>0\},$ $\|x^*-\bar{x}\|\geq \bar{r}$. By $x^*\in\overline{\mathcal{B}(\bar{x}, \bar{r})}$, $\|x^*-\bar{x}\|\leq \bar{r}$. Hence, \begin{equation}\label{27uuu}
	\|x^*-\bar{x}\|^2 = \bar{r}^2=\inf_{\{g=0\}}\|x-\bar x\|^2.
	\end{equation}
Due to the assumption $\{g=0\}=\{h=0\}$, \eqref{27uuu} applies likewise to $\{h=0\}.$
	
	Next, we claim that $\nabla g(x^*)\ne 0$ and $\nabla h(x^*)\ne 0.$
	
	Suppose contrarily that there is some $x^*\in \overline{\mathcal{B}(\bar{x}, \bar{r})}\cap \{g=0\}$ such that $\nabla g(x^*)= 0.$
	Then, with a suitable change of variables, $g(x)$ can only be of the canonical form \eqref{form:1}; or \eqref{form:2}; or \eqref{form:4} with $\eta=0$. That is, $$g(x)=\bar{\delta}(-x_1^2-\cdots -x_k^2)+\delta(x_{k+1}^2+\cdots +x_m^2)+\theta,$$ where $\bar{\delta}, {\delta}\in \{0, 1\}, \theta\in \mathbb{R}.$ Moreover, since $\nabla g(x^*)= 0$, one has
	\begin{equation*}\label{}
	x^*=(x^*_1, \cdots, x^*_k, x^*_{k+1}, \cdots, x^*_m, \cdots, x^*_n )^T \mbox{ where } \left\{\begin{array}{ll}
	x^*_1=\cdots =x^*_k=0,~\mbox{ if } \bar{\delta}=1; \\
	x^*_{k+1}=\cdots =x^*_m=0, ~\mbox{ if } {\delta}=1.
	\end{array}\right.\end{equation*}
	Due to $g(x^*)=0$, $\theta$ must be zero.
	Assumption 4 further implies that $\bar{\delta}$ and $\delta$ must be $1$ so that
	$x^*=(0,\cdots, 0, x^*_{m+1}, \cdots, x^*_n)^T$
	and $$g(x)=(-x_1^2-\cdots -x_k^2)+(x_{k+1}^2+\cdots +x_m^2).$$
	
	Now, let
	\begin{equation*}\label{}
	\bar{y}=(\bar{y}_1,\cdots,\bar{y}_m,\bar{y}_{m+1},\cdots,\bar{y}_n )^T \mbox{ with } \left\{\begin{array}{ll}
	\bar{y}_i=0,~\mbox{ for } i\in \{1, \cdots, m\}; \\
	\bar{y}_i=\bar{x}_i, ~\mbox{ for } i\in \{m+1, \cdots, n\},
	\end{array}\right.\end{equation*}
	where $\bar{x}_i$ is the coordinates of $\bar{x}$.
	Then $\bar{y}\in \{g=0\}$ and $\|\bar{y}-\bar{x}\|=\sqrt{\bar{x}_1^2+\cdots+\bar{x}_m^2}$.
	However,  $$\|x^*-\bar{x}\|=\sqrt{\bar{x}_1^2+\cdots+\bar{x}_m^2+\sum_{i=m+1}^{n}(x^*_{i}-\bar{x}_{i})^2}\geq \|\bar{y}-\bar{x}\|,$$ which together with \eqref{27uuu} implies that $x^*_i=\bar{x}_i$ for all $i\in \{m+1, \cdots, n\}$. In other words,
	\begin{equation}\label{}
	x^*=(0,\cdots, 0, \bar{x}_{m+1}, \cdots, \bar{x}_n)^T.
	\end{equation}
	By $g(\bar{x})>0$, there exists some $j\in\{k+1,\cdots, m\}$ such that $\bar{x}_j\ne 0.$ Without loss of generality, let $j=m$ and
	\begin{equation}\label{29sss}
	\bar{x}_m\ne 0.
	\end{equation}
	
	Consider
	\begin{equation*}
	\hat{x}=(\hat{x}_1,\cdots,\hat{x}_m,\hat{x}_{m+1},\cdots,\hat{x}_n )^T \mbox{ with } \left\{\begin{array}{ll}
	\hat{x}_i=\dfrac{\bar{x}_1+\bar{x}_m}{2},~\mbox{ for } i=1,m; \\
	\hat{x}_i=0, ~\mbox{ for } i\in \{2, \cdots, m-1\};\\
	\hat{x}_i={\bar x}_i, ~\mbox{ for } i\in \{m+1, \cdots, n\}.
	\end{array}\right.\end{equation*}
	It is easy to see that $\hat{x}\in \{g=0\}$ and
	\begin{equation*}\label{}	\|\bar{x}-\hat{x}\|=\sqrt{\big(\dfrac{\bar{x}_1-\bar{x}_m}{2}\big)^2
		+\big(\dfrac{\bar{x}_m-\bar{x}_1}{2}\big)^2+\sum_{i=2}^{m-1}\bar{x}_i^2}.
	\end{equation*}
	On the other hand,
	\begin{equation*} 
	\|\bar{x}-x^*\|=\sqrt{\bar{x}_1^2+\bar{x}_m^2+\sum_{i=2}^{m-1}\bar{x}_i^2}.
	\end{equation*}
	Obviously, $\|\bar{x}-x^*\|\geq \|\bar{x}-\hat{x}\|$ and $\|\bar{x}-x^*\|= \|\bar{x}-\hat{x}\|$ if and only if $\bar{x}_1=-\bar{x}_m$.
	\begin{itemize}
		\item If $\bar{x}_1\ne-\bar{x}_m,$ then $\|\bar{x}-x^*\|> \|\bar{x}-\hat{x}\|$, contradicting to \eqref{27uuu}.		
		\item Otherwise, $\bar{x}_1=-\bar{x}_m.$ In this case, choose
		\begin{equation*} 
		\check{x}=(\check{x}_1,\cdots,\check{x}_m,\check{x}_{m+1},\cdots,\check{x}_n )^T \mbox{ with } \left\{\begin{array}{ll}
		\check{x}_1=-\check{x}_m=\dfrac{\bar{x}_1-\bar{x}_m}{2}; \\
		\check{x}_i=0, ~\mbox{ for } i\in \{2, \cdots, m-1\};\\
		\check{x}_i={\bar x}_i, ~\mbox{ for } i\in \{m+1, \cdots, n\}.
		\end{array}\right.\end{equation*}
		Then, $\check{x}\in \{g=0\}$ and
		\begin{eqnarray}\label{l1} \|\bar{x}-\check{x}\|&=&\sqrt{\big(\dfrac{\bar{x}_1+\bar{x}_m}{2}\big)^2+
			\big(\dfrac{\bar{x}_m+\bar{x}_1}{2}\big)^2+\sum_{i=2}^{m-1}\bar{x}_i^2} \nonumber\\
		&=&\sqrt{\sum_{i=2}^{m-1}\bar{x}_i^2} \nonumber\\
		&<&\sqrt{\bar{x}_1^2+\bar{x}_m^2+\sum_{i=2}^{m-1}\bar{x}_i^2}=\|\bar{x}-x^*\|, \label{212uuu}
		\end{eqnarray}
		where the inequality \eqref{212uuu} follows from \eqref{29sss}. So we have proved that $\|\bar{x}-\check{x}\|<\|\bar{x}-x^*\|$, which again contradicts to \eqref{27uuu} since $\check{x}\in \{g=0\}$.
	\end{itemize}
	
	In summary, we have proved that $\nabla g(x^*)\not=0,$ which indicates that  $$x^*\in\arg\min\{\|x-\bar{x}\|^2: g(x)=0\}$$ is a regular point of $g(x)=0.$ By the KKT condition, there exists a $\rho\in\mathbb{R}$ such that
	$$2(x^*-\bar x)+\rho\nabla g(x^*)=0.$$
	Since $x^*\ne\bar x,$ $\rho\not=0.$ Then, $\nabla g(x^*)=t^*\overrightarrow{x^*\bar{x}}$
	with $t^*=\frac{2}{\rho}.$ Since $g(x^*)=0$ and $~g(\bar x)>0,$  there must be $t^*=\frac{2}{\rho}>0.$
	By a similar argument, we can show the same for $\nabla h(x^*).$ The proof is thus complete.
	\hfill$\Box$
	\vskip0.3cm
	
	
\noindent {\bf Proof for Theorem \ref{lm31aa}: the unsolvability of \eqref{Unsolvable-system} under Assumptions 1 - 5.}
	\vskip0.3cm
\proof
	According to Assumption 2, the following two cases must both happen.
	\begin{itemize}
		\item[($\sharp$)] either $\{g=0\}\subset \{h\leq 0\}$ or $\{g=0\}\subset \{h\geq 0\},$ in which case, by Lemma \ref{lemma3-3}, $\{h=0\}$ separates neither $\{g<0\}$ nor $\{g>0\}.$
		\item[($\flat$)] either $\{h=0\}\subset \{g\leq 0\}$ or $\{h=0\}\subset \{g\geq 0\},$ in which case, by Lemma \ref{lemma3-3}, $\{g=0\}$ separates neither $\{h<0\}$ nor $\{h>0\}.$
	\end{itemize}
	
	In ($\flat$), when $\{h=0\}\subset \{g\geq 0\}$ happens, either case in ($\sharp$) implies that $\{h=0\}$ cannot separate $\{g<0\}.$ By Lemma \ref{th23sss-1}, there is a
	$\lambda_1\in \mathbb{R}$ such that
	\begin{equation}
	g(x)+\lambda_1 h(x)\geq 0. \label{z33hjk}
	\end{equation}
	Analogously, when $\{h=0\}\subset \{-g\geq 0\}$ happens in ($\flat$),
	($\sharp$) implies that $\{h=0\}$ cannot separate $\{-g<0\}.$ By Lemma \ref{th23sss-1} again, there exists $\lambda_2\in \mathbb{R}$ such that
	\begin{equation}
	-g(x)+\lambda_2 h(x)\geq 0. \label{z34hjk}
	\end{equation}
	
	Under the Slater condition (Assumption 4), we can further apply the classical S-lemma to see that there
	are 4 possibilities:
	
	- \eqref{z33hjk} implies that
	\begin{equation}
	\{h\leq 0\}\subset \{g\geq 0\} \hbox{ when }\lambda_1\geq 0; \label{z37hjk}
	\end{equation}
	or
	\begin{equation}
	\{h\geq 0\}\subset \{g\geq 0\} \hbox{ when }\lambda_1< 0; \label{z38hjk}
	\end{equation}
	
	- \eqref{z34hjk} implies that
	
	\begin{equation}
	\{h\leq 0\}\subset \{g\leq 0\} \hbox{ when }\lambda_2\geq 0; \label{z35hjk}
	\end{equation}
	or
	\begin{equation}
	\{h\geq 0\}\subset \{g\leq 0\} \hbox{ when }\lambda_2< 0. \label{z36hjk}
	\end{equation}

	\noindent Case 1:  \eqref{z35hjk} cannot occur since it implies that $\mathcal{D}=\{h\leq 0\}$, which contradicts to Assumption $3$.
	
	\noindent Case 2:  By the same argument, \eqref{z38hjk} cannot occur, either, since
	it implies that $\{g< 0\}\subset \{h < 0\}$ so that $\{g\leq 0\}\subset \{h \leq 0\}$ and thus
	$\mathcal{D}=\{g\leq 0\}.$
	
	\noindent Case 3:  In case of \eqref{z36hjk}, there is $\{g> 0\}\subset\{h< 0\}$. Then, the system
	$$  \begin{cases}
	\begin{array}{l}
	g(x)>0\\
	h(x)\geq 0
	\end{array}\end{cases}$$ has no solution.
	
	\noindent Case 4:  In case of \eqref{z37hjk}, we have
	$$\mathcal{D}=\{h \leq 0\}\cap \{g \leq 0\}\subset \{g\geq 0\}\cap \{g \leq 0\}=\{g=0\}.$$ By Assumption 1, $\mathcal{D}=\{g=0\}.$
	On the other hand, \eqref{z37hjk} also implies that $\{g\leq 0\}\subset \{h \geq 0\}$. So $\mathcal{D}=\{h \leq 0\}\cap \{g \leq 0\}\subset \{h\leq 0\}\cap \{h \geq 0\}=\{h=0\}.$ By Assumption 1, $\mathcal{D}=\{h=0\}.$ That is, when the case \eqref{z37hjk} happens,
	\begin{equation}\label{213ggg}
	\mathcal{D}=\{g=0\}=\{h=0\}\not=\emptyset~~ (\hbox{by Assumption 3}).
	\end{equation}
	
	Our goal is to prove that, under Case 4, the system of two quadratic inequalities
	$$ \begin{cases}
	\begin{array}{l}
	g(x)>0\\
	h(x)\geq 0
	\end{array}\end{cases}$$ has no solution.
	Suppose on the contrary that this system has a solution $\bar{x}$ such that $g(\bar{x})>0, h(\bar{x})\ge0.$ Since $g(\bar{x})\not=0,$ by \eqref{213ggg}, $h(\bar{x})\not=0.$ Namely, one has
	\begin{equation*}
	g(\bar{x})>0,~ h(\bar{x})>0.
	\end{equation*}
	By Lemma \ref{l24dd}, there exists a largest ball $\mathcal{B}(\bar{x}, \bar{r}),$ centered at
	$\bar{x}$ with $\bar{r}<\infty$ such that
	\begin{equation}\label{214aaa}
	\mathcal{B}(\bar{x}, \bar{r})\subset \{g>0\}\cap \{h>0\}.
	\end{equation}
	and there exists some point $x^*\in \overline{\mathcal{B}(\bar{x}, \bar{r})}$ such that
	\begin{equation*}
	g(x^*)=h(x^*)=0, \nabla g(x^*)=t^*d\ne 0 \text{ and } \nabla h(x^*)=s^*d\ne 0
	\end{equation*}
	where $d=\overrightarrow{x^*\bar{x}}.$
	Then, by Taylor's expansion,
	\begin{eqnarray}
	g(x^*+\alpha \nabla g(x^*)) &=& g(x^*)+ \alpha\nabla g(x^*)^T\nabla g(x^*)
	+\alpha^2\big(\nabla g(x^*)^TB\nabla g(x^*)\big)\nonumber\\
	&=& \alpha \|t^*d\|^2+\alpha^2(t^*)^2\big(d^TBd\big)\label{Taylor-g};
	\end{eqnarray}
	and
	\begin{eqnarray}
	h(x^*+\beta \nabla h(x^*)) &=& h(x^*)+ \beta\nabla h(x^*)^T\nabla h(x^*)
	+\beta^2\big(\nabla h(x^*)^TC\nabla h(x^*)\big)\nonumber\\
	&=& \beta \|s^*d\|^2+\beta^2(s^*)^2\big(d^TCd\big)\label{Taylor-h}.
	\end{eqnarray}
	From \eqref{Taylor-g} and \eqref{Taylor-h},
	for all sufficiently small negative values of $\alpha,\beta\in(\gamma,0),~\gamma<0,$ one has
	$$g(x^*+\alpha \nabla g(x^*))<0;~~h(x^*+\beta \nabla h(x^*))<0.$$
	Let us assume that $t^*<s^*$ and take $\alpha=\frac{\gamma}{2},~\beta=\frac{\gamma t^*}{2s^*},~\ddot{x}=x^*+\frac{\gamma t^*}{2}.$ Then, $$g(\ddot{x})<0,~~ h(\ddot{x})<0,$$
	so that $\ddot{x}\in\mathcal{D}.$ However, it
	contradicts to \eqref{213ggg} that $\mathcal{D}=\{g=0\}=\{h=0\}.$
	
	In summary, based on ($\flat$), we have shown that the system $\{g(x)>0, h(x) \geq 0\}$ has no solution. The same analysis to base on ($\sharp$) will lead to the fact that the system
	$\{g(x)\ge0, h(x) >0\}$ cannot have a solution, either.  \hfill$\Box$

	\section{Solving $\rm(P)$ under Assumptions 1 and 2}\label{sec:results}
	
	Denote $\nu^*(\hbox{[Non-Alter]})$ to be the optimal value of a problem (P) satisfying Assumptions 1 and 2. The key question in this section is to investigate the feasibility problem:
	
	Given any real number $\gamma\in\mathbb{R},$ can we check efficiently whether the sublevel set of the objective function $\{f<\gamma\}$ is feasible or not on the constraint set $\mathcal{D}?$ That is, determine whether or not $\{f<\gamma\}\cap \mathcal{D}=\emptyset.$
	
	The feasibility problem stated above links to the optimal value of a problem in [Non-Alter] by
	$$\nu^*(\hbox{[Non-Alter]})=\sup\{\gamma:~\{f<\gamma\}\cap \mathcal{D}=\emptyset\}.$$
	
	For an infeasible solution $x$ to $\mathcal{D},~ x\not\in\mathcal{D},$ either $g(x)>0,$ or $h(x)>0,$ or both happen. If, in addition, we have Assumptions 3 - 5, by the unsolvability Theorem \ref{lm31aa}, when $g(x)>0,$ $h(x)$ must be negative, and vice verse. In other words,
	\begin{equation}\label{z311ghjk}
	\{f<\gamma\}\cap \mathcal{D}=\emptyset \Longrightarrow	\{f<\gamma\}\subset \big(\{g<0\}\cap\{h> 0\}\big)\cup\big(\{g> 0\}\cap\{h< 0\}\big).
	\end{equation}
	
	The following theorem shows that, in case of \eqref{z311ghjk}, $\{f<\gamma\}$ must lie, either entirely in $\big(\{g<0\}\cap\{h> 0\}\big)$, or entirely on the other side $\big(\{g> 0\}\cap\{h< 0\}\big),$  except for a case when both $g,~h$ are affine. More importantly, which side $\{f<\gamma\}$ would reside has nothing to do with the choice of $\gamma.$ It is solely determined by the characteristic of the objective function $f.$
	
	\begin{theorem}\label{l32rtyu} Suppose that neither $g$ nor $h$ is affine and Assumptions $1,2,3,4$ and $5$ hold true. Let $\Gamma=\big\{\gamma\in \mathbb{R}: \{f<\gamma\}\cap \mathcal{D}=\emptyset\big\}.$ Then,
		\begin{eqnarray}
		\text{either }\	(\forall \gamma\in\Gamma),~\{f<\gamma\}\subset \{g< 0\}\cap \{h> 0\}; \label{233poiu}\\
		\text{ or }\ (\forall \gamma\in\Gamma),~	\{f<\gamma\}\subset \{g> 0\}\cap \{h< 0\}. \label{234poiu}
		\end{eqnarray}
	\end{theorem}
	
	\proof
	
	Let $\gamma\in \Gamma.$ Since $f$ is quadratic, the sublevel set $\{f<\gamma\}$ can have at most two connected components. For detailed proof, please refer to Lemma 1 in Quang and Sheu \cite[2018]{Quang-Sheu18}). We first observe that each connected component $L$ of $\{f<\gamma\}$ must lie either entirely in $\{g<0\}$ or entirely in $\{g>0\}$. Otherwise, by Intermediate value theorem, there exists $x^*\in L\subset \{f<\gamma\}$ such that $g(x^*)=0.$ By Lemma \ref{lm31aa}, $\{g\geq 0\}\cap\{h> 0\}=\emptyset,$ so that $h(x^*)\le0$ and thus
	$x^*\in \mathcal{D}.$ It is impossible as $x^*\in\{f<\gamma\}\cap \mathcal{D}=\emptyset.$ From \eqref{z311ghjk}, we conclude that if $L$ is a connected component of $\{f<\gamma\}$,
	either $L\subset \big(\{g<0\}\cap\{h> 0\}\big)$ or $L\subset\big(\{g> 0\}\cap\{h< 0\}\big).$
	
	Now suppose that $\{f<\gamma\}$ consists of two connected components, $L_1$ and $L_2,$ such that
	$$L_1\subset\big(\{g<0\}\cap\{h> 0\}\big) \hbox{ and } L_2\subset\big(\{g> 0\}\cap\{h< 0\}\big).$$
	Then, $g(L_1)g(L_2)<0$ and $h(L_1)h(L_2)<0,$ indicating that
	both $\{g=0\}$ and $\{h=0\}$ separate $\{f<\gamma\}$. By \cite[Theorem 1]{Quang-Sheu18}, both $\{g=0\}$ and $\{h=0\}$ are affine functions, which contradicts to the assumption of the theorem.
	So, we have proved that, $\forall \gamma\in \Gamma,$
	\begin{equation}\label{z237kjhg}
	\text{either }\,\{f<\gamma\} \subset \big(\{g >0\}\cap\{h<0\}\big);
	\text{ or }\,\{f<\gamma\} \subset\big(\{g <0\}\cap\{h>0\}\big).\end{equation}
	
	Finally, we have to show, if $\emptyset\ne \{f<\gamma^*\} \subset \big(\{g >0\}\cap\{h<0\}\big)$ for some $\gamma^*\in \Gamma,$ then
	\begin{equation}\label{z312dsa}
	\{f<\gamma\}\subset (\{g >0\}\cap\{h<0\}), ~~ \forall \gamma \in \Gamma.
	\end{equation} 
	
	- It is obvious that $\{f<\gamma\}\subset\{f<\gamma^*\}$ if $\gamma<\gamma^*.$ Therefore, \eqref{z312dsa} holds for all $\gamma \in \Gamma, \gamma\leq \gamma^*$.
	
	- If $\gamma>\gamma^*$, then $\{f<\gamma\}\supset\{f<\gamma^*\}$. Suppose contrarily that
	$\{f<\gamma\}\subset \big(\{g <0\}\cap\{h>0\}\big).$ It would imply that $\{f<\gamma^*\}\subset (\{g <0\}\cap\{h>0\}),$ too. This contradicts to $\emptyset\ne \{f<\gamma^*\} \subset \big(\{g >0\}\cap\{h<0\}\big)$. By \eqref{z237kjhg}, it follows that \eqref{z312dsa} holds for all $\gamma \in \Gamma, \gamma> \gamma^*$.
	The proof is thus completed.
	\hfill$\Box$
	
	With Theorem \ref{l32rtyu}, we have the following $\mathcal{S}$-procedure for two quadratic inequalities.

	\begin{theorem}\label{l33rtyu} Suppose that Assumptions $1,2,3,4$ and $5$ hold and
		neither $g$ nor $h$ is affine. The following two statements are equivalent for any given $\gamma\in\mathbb{R}$.
		
		${\rm(S_1)}$ $(\forall x\in \mathbb{ R }^n)~\big(g(x)\le0,~h(x)\le0\big)~\Longrightarrow~f(x)\geq \gamma$.
		
		${\rm(S_2)}$ $(\exists \lambda_1\ge0, \lambda_2\ge0)$ such that $f(x)-\gamma+\lambda_1g(x)+\lambda_2h(x)\geq 0,~~ \forall x\in \mathbb{ R }^n.$	
	\end{theorem}
	
	\proof It is trivial that ${\rm(S_2)}$ implies ${\rm(S_1)}$. We just need to prove ${\rm(S_1)}$ implies ${\rm(S_2)}$.  When we assume ${\rm(S_1)}$, by Theorem \ref{l32rtyu}, one of \eqref{233poiu} and \eqref{234poiu} must occur.
	
	- In case of \eqref{233poiu}: $\{f<\gamma\}\subset\{h>0\}.$ Equivalently, $\{h\le0\}\subset\{f-\gamma\geq 0\}$. By S-lemma, there exists
	$\lambda\geq 0$ such that $f(x)-\gamma+\lambda h(x)\geq 0, \forall x\in \mathbb{ R }^n.$ Choose $\lambda_1=0, \lambda_2=\lambda$ for ${\rm(S_2)}$ to follow.
	
	- In case of \eqref{234poiu}: the same argument as in \eqref{233poiu} applies.
	\hfill$\Box$
	
	\begin{theorem}\label{thm34hjkl} Suppose that Assumptions $1, 2$ hold true. Then, $\nu^*(\hbox{[Non-Alter]})$ can be computed in polynomial time.	
	\end{theorem}
	
	\proof We divide the proof into the following cases:
	
	$\odot$ If Assumption 3 is violated, then either $\mathcal{D}=\emptyset$
	in which case [Non-Alter] is infeasible; or
	$\mathcal{D}=\{g\leq 0\}$ or $\mathcal{D}=\{h\leq 0\},$ in which case [Non-Alter] is reduced to [QP1QC].
	
	
	$\odot$	Assumption 4 is the two-sided Slater condition, which
	is {\it not} essential in proving strong duality of [Non-Alter]. Suppose $g(x)$ fails the Slater condition. Then, $g(x)\geq 0,~\forall x\in \mathbb{R}^n;$ $g(x)$ is convex, and $$\{g\leq 0\}=\{g=0\}=\{x\in\mathbb{R}^n|~Bx+b=0\}$$ is a hyperplane in $\mathbb{R}^n.$
	In this case, $\mathcal{D}$ is the intersection of $\{h\leq 0\}$ with a hyperplane, which reduces [Non-Alter] to [QP1QC] over $\mathbb{R}^{n-1}.$ Similar arguments can be applied to $-g,~h,~-h$ when any of them violates the Slater condition.
	
	$\odot$	Assumption 5, once failed, means that $g(x)=-x_1^2+1,~h(x)={\bar{c}_1}x_1+{\bar{c}_0},~\bar{c}_1>0,~\bar{c}_0=\pm \bar{c}_1\not=0$ after coordinate transformation.
	
	- If $\bar{c}_0=\bar{c}_1,$  $\mathcal{D}=\{x\in\mathbb{R}^n:~x_1 \le-1\}.$ Then, [Non-Alter] is a [QP1QC].
	
	- If $\bar{c}_0=-\bar{c}_1,$  $\mathcal{D}=\{x\in\mathbb{R}^n:~x_1 \le-1\}\cup \{x\in\mathbb{R}^n:~x_1=1\}.$ In this case
	$$\nu^*(\hbox{[Non-Alter]})=\min\big\{\min\limits_{x_1\le -1}f(x);
	\min\limits_{x_1=1}f(x) \big\},$$
	which can be solved in polynomial time.
	
	$\odot$ When both $g,~h$ are affine and satisfy Assumptions 1 and 2,
	we may set $g(x)=2b^Tx+b_0, h(x)=2c^Tx+c_0$. By Assumption 2, we have $c=tb, t\ne 0$.
	
	- If $t>0,$ then $\mathcal{D}=\{g\leq 0\}$ or $\mathcal{D}=\{h\leq 0\}.$ [Non-Alter] becomes [QP1QC].
	
	- If $t<0,$ then $\mathcal{D}=\{x: \frac{-c_0}{2t}\leq b^Tx \leq -\frac{b_0}{2}\}=\{x: (b^Tx+\frac{c_0}{2t})(b^Tx+\frac{b_0}{2})\leq 0\}$ so that [Non-Alter] becomes [QP1QC] again.

	$\odot$	Finally, we arrive the latest case when Assumptions 1 - 5 hold, and neither $g$ nor $h$ is affine so that the $\mathcal{S}$-procedure Theorem \ref{l33rtyu} applies. Then (P) and its SDP relaxation has no gap. Indeed, with the standard notation defined in \eqref{M-notation}, we have
	\begin{eqnarray} &&\nu^*(\hbox{[Non-Alter]})\\
	&=&\inf\left\{f(x) :~ x\in\mathcal{D}\right\}\nonumber \\
	&=&\sup\left\{\gamma :~ \{f<\gamma\}\cap\mathcal{D}=\emptyset\right\} ~\label{32qwer--1} \\
	&=&\sup_{(\lambda_1,\lambda_2)\geq 0}\left\{\gamma :~ f(x)-\gamma+\lambda_1g(x)+\lambda_2h(x)\geq 0, \forall x\in\mathbb{ R }^n \right\}  ~ \label{32qwer}\\
	&=&\sup_{(\lambda_1,\lambda_2)\geq 0}\left\{\gamma :Z=M(f)-\gamma \begin{bmatrix}1& \bar{0}^T \\ \bar{0} & [0] \end{bmatrix}+\lambda_1 M(g)+\lambda_2 M(h),Z\succeq 0 \right\}\label{33qwer}\\
	&\leq&\min\left\{M(f)\bullet X :~ M(g)\bullet X \leq 0,  M(h)\bullet X \leq 0, X_{00}=1, X\succeq 0 \right\}\label{34qwer}\\
	&\leq&\min\left\{M(f)\bullet X :~ M(g)\bullet X \leq 0,  M(h)\bullet X \leq 0, X=\begin{bmatrix}1 \\x \end{bmatrix}\begin{bmatrix}1 \\x \end{bmatrix}^T \right\}\label{35qwer}\\
	&=&\min\left\{f(x) :~ g(x) \leq 0,  h(x) \leq 0\right\} \label{36qwer}\\
	&=&\nu^*(\hbox{[Non-Alter]}) \nonumber.
	\end{eqnarray}
	
	Note that \eqref{32qwer--1} is the statement ${\rm(S_1)}$ in Theorem \ref{l33rtyu};
	\eqref{32qwer} is the statement ${\rm(S_2)}$ in Theorem \ref{l33rtyu};
	\eqref{33qwer} is the SDP reformulation of \eqref{32qwer}; \eqref{34qwer} follows from the conic weak duality; \eqref{35qwer} follows from the fact that
	$$ X=\begin{bmatrix}1 \\x \end{bmatrix}\begin{bmatrix}1 \\x \end{bmatrix}^T\succeq0 \mbox{ with } X_{00}=1;$$
	and \eqref{36qwer} follows from $$M(f)\bullet \begin{bmatrix}1 \\x \end{bmatrix}\begin{bmatrix}1 \\x \end{bmatrix}^T=\begin{bmatrix}1 \\x \end{bmatrix}^TM(f)\begin{bmatrix}1 \\x \end{bmatrix}=f(x).$$
	Therefore, the optimal value $\gamma^*$ of (P) can be solved by solving the SPD \eqref{33qwer}. All the procedures discussed above are polynomial.
	The theorem is thus proved. \hfill$\Box$
	
	\subsection{Finding an optimal solution for [Non-Alter]}	
	
	As for the optimal solution of [Non-Alter], unlike [YZ] in which the primal-dual Slater conditions were assumed so that the rank-one decomposition can be applied for finding an optimal solution from the solution of the SDP relaxation (SP), here we have to go via another approach.
	
	In general, the objective function in [Non-Alter] can be unbounded from below on the set of feasible solutions, or bounded from below but not attainable. Here we aim to find an optimal solution $x^*$ for [Non-Alter] when the optimal value $\nu^*(\hbox{[Non-Alter]})$ is finite and attainable.
	
	We divide into two cases.
	
	\noindent $\blacktriangleright$ ~$-\infty<\nu^*(\hbox{[Non-Alter]})=\inf \{f(x): x\in\mathbb{ R }^n\}$:
	
	When the unconstrained minimization of $f$ over $\mathbb{R}^n$ happens to equal
	$\nu^*(\hbox{[Non-Alter]}),$ the objective function $f$ must be convex. An optimal solution $x^*$ to [Non-Alter] can be found by solving
	a solution from the following system:
	\begin{equation}\label{3.7poiu}
	\left\{\begin{array}{ll}
	f(x)=\nu^*(\hbox{[Non-Alter]}),  \\
	g(x)\leq 0,  \\
	h(x)\leq 0.
	\end{array}\right.\end{equation}
	
	The first equation in \eqref{3.7poiu} consists of all the minimizers of $f(x)$ on $\mathbb{ R }^n$ and this set is equivalent to
	\begin{equation}
	\left\{x: \min _{x\in\mathbb{R}^n} f(x)=\nu^*(\hbox{[Non-Alter]})\right\}=\left\{-A^{+} a+Z y: y \in \mathbb{R}^{m}\right\}, \nonumber
	\end{equation}
	where $A^+$ is the pseudo-inverse of $A$ and $Z \in \mathbb{ R }^{n\times m}$ is a matrix basis of $\mathcal{N}(A)$, assuming the null space of $A$ has rank $m.$
	
	The system \eqref{3.7poiu} now reduces to finding an $y\in\mathbb{R}^m$ satisfying
	\begin{equation}\label{3.7poiu-1}
	\left\{\begin{array}{ll}
	\bar{g}(y)=g(-A^+a+Zy)\leq 0,  \\
	\bar{h}(y)=h(-A^+a+Zy)\leq 0.
	\end{array}\right.\end{equation}
	in which both $\bar{g}(y), \bar{h}(y)$ are quadratic in $y.$ By solving the following [QP1QC]
	\begin{equation}\label{1uu}
	\inf\{\bar{g}(y): \bar{h}(y)\leq 0\},
	\end{equation}
	we can get an optimal solution $y^*$ to \eqref{1uu}. Then, $x^*=-A^{+} a+Z y^*$ is optimal to [Non-Alter].
	
	\noindent $\blacktriangleright$ ~$\nu^*(\hbox{[Non-Alter]})>\inf \{f(x): x\in\mathbb{ R }^n\}=\gamma^*$:
	
	We first notice that $\{f<\nu^*(\hbox{[Non-Alter]})\}\ne \emptyset$ and
	$\{f<\nu^*(\hbox{[Non-Alter]})\}\cap \mathcal{D}=\emptyset.$ According to Theorem \ref{l32rtyu},
	one of the following two cases
	\begin{eqnarray} 
	\{f<\nu^*(\hbox{[Non-Alter]})\}\subset \{g< 0\}\cap \{h> 0\}, \label{3.9poiu} \\
	\text{ or }~ \{f< \nu^*(\hbox{[Non-Alter]})\}\subset \{g> 0\}\cap \{h< 0\}. \label{3.10poiu}
	\end{eqnarray}
	occurs. This is easy to determine by taking a point $x^0\in\{f<\nu^*(\hbox{[Non-Alter]})\}$ and checking the values of $g(x)$ and $h(x)$ at $x^0$. Without loss of generality, let us assume \eqref{3.9poiu}. Then,
	$$\{f<\nu^*(\hbox{[Non-Alter]})\}\cap \{h\le0\}=\emptyset.$$
	It implies that
	\begin{equation}\label{2-oo}
	f(x^*)=\nu^*(\hbox{[Non-Alter]})\le \inf\{f(x):~h(x)\le0\}.
	\end{equation}
	Obviously, $\inf\{f(x):~h(x)\le0\}\le\nu^*(\hbox{[Non-Alter]}),$ which, together with
	\eqref{2-oo} gives
	\begin{equation}\label{3-oo}
	f(x^*)=\nu^*(\hbox{[Non-Alter]})=\inf\{f(x):~h(x)\le0\}.
	\end{equation}
	As $x^*\in\mathcal{D},$ $x^*$ is an optimal solution to $\inf\{f(x):~h(x)\le0\}$ and $\inf\{f(x):~h(x)\le0\}$ is attainable. Let us solve the [QP1QC] problem $\inf\{f(x):~h(x)\le0\}$ and assume the point $x_h$ to be the optimal solution obtained.
	
	We claim that $h(x_h)\not<0.$ Otherwise, $x_h$ is an interior point so that
	\begin{equation}\label{4-oo}
	f(x_h)=\inf\{f(x):~h(x)\le0\}=\inf\{f(x): x\in \mathbb{R}^n\}.
	\end{equation}
	By \eqref{3-oo} and \eqref{4-oo}
	$$f(x_h)=\nu^*(\hbox{[Non-Alter]})=\gamma^*,$$
	which contradicts to the assumption at the beginning of this case. So, $h(x_h)=0.$
	However, by the unsolvability Theorem \ref{lm31aa}, $g(x_h)\not>0.$ There must be
	$g(x_h)\le0.$ In other words, $x_h\in\mathcal{D}.$ The point $x_h$ is optimal to [Non-Alter].

	\section{Final Discussions}
	
	Much of the quadratically constrained quadratic programming literature assumes that both (SP) and (SD) satisfy the Slater condition and that $g(x)$ and $h(x)$ do not share any common root. Our paper removes those obstacles by exploring topological relations (non-alternativeness) between the (sup-/sub-) level sets of $g$ and $h.$ We have shown that [Non-Alter] covers [QP1QC],~[QP1EQC], [GTRS], and [YZ], but not [CDT].
	
	Our final comment goes to [HQPD]. It is easy to find an example in [HQPD] fails Assumption 1. For example, let $f(x)=x_1^2-x_2^2$, $g(x)=x_1^2+x_2^2-1, h(x)=-x_1^2+1$, then $A+2B\succ 0, 2B+C\succ 0.$ This example belongs to [HQPD]. However,
	$\mathcal{D}=\{g\leq 0\}\cap \{h\leq 0\}=\{(-1,0)^T, (1, 0)^T\}\subset \{g=0\}$ but $\mathcal{D}\ne \{g=0\}$ and  $\mathcal{D}=\{g\leq 0\}\cap \{h\leq 0\}=\{(-1,0)^T, (1, 0)^T\}\subset \{h=0\}$ but $\mathcal{D}\ne \{h=0\}$ so they do not satisfy Assumption 1. On the other hand, $\{h=0\}\subset\{g\ge0\}$ and $\{g=0\}\subset\{h\ge0\}.$ Assumption 2 is satisfied by this example.
	
	Now consider another example. Let $f(x)=x_1^2+x_2^2,~g(x)=2x_1^2+x_2^2-9, h(x)=x_1^2+2x_2^2-9$. This example belongs to [HQPD] because $A,B,C\succ0.$
	It is easy to see that $\mathcal{D}\ne \{g=0\}$, $\mathcal{D}\ne \{h=0\}$ so they satisfy Assumption 1. However,
	$\{g=0\}\not\subset \{h\le0\}$, $\{g=0\}\not\subset \{h\ge0\}$, $\{h=0\}\not\subset \{g\le0\}$, $\{h=0\}\not\subset \{g\ge0\}$ so they do not satisfy Assumption 2.
	
	Interestingly, with all $f,g,h$ in homogeneous forms having a pd pencil, we could {\it not} find an example which fails both Assumption 1 and Assumption 2. Polyak's result \cite[1998]{Polyak98} is a direct consequence of Brickman \cite[1961]{Brickman}. We wish that we shall be able to extend their results for more general [HQPD] in the future.

	\section*{Acknowledgements}
	
	Huu-Quang, Nguyen's research work was sponsored partially by Taiwan Ministry of Science and Technology grant number: MOST 108-2811-M-006-537 \\
	Ruey-Lin Sheu's research work
	was sponsored partially by Taiwan Ministry of Science and Technology grant number: MOST 107-2115-M-006-011-MY2.

\end{document}